%% file: main.tex
\documentclass[12pt]{article}

\usepackage[a4paper,
	left=2.5cm,
	right=2.5cm,
	top=3cm,
	bottom=4cm,
	footskip=1.5cm,
	twoside,
]{geometry}

\input{preamble.tex}

\begin{document}
\baselineskip=1.5em
\begin{center}
	{\Large {\bf Fifth-order finite volume derivative-based Hermite WENO scheme with unified stencils for hyperbolic conservation laws}}
\end{center}

\begin{center}
	Peiwen Chen\footnote{School of Mathematical Sciences,  Xiamen University, Xiamen, Fujian 361005, P.R. China. E-mail: chenpw@stu.xmu.edu.cn.},
	and
	Zhuang Zhao\footnote{Corresponding author. School of Mathematical Sciences and Fujian Provincial Key Laboratory of Mathematical Modeling and High-Performance Scientific Computing, Xiamen University, Xiamen, Fujian 361005, P.R. China. E-mail: zzhao@xmu.edu.cn.}
\end{center}
\input{./sections/sec-abstract.tex}

\input{./sections/sec-introduction.tex}

\input{./sections/sec-description.tex}

\input{./sections/sec-numerical-tests.tex}

\input{./sections/sec-conclusion.tex}

\input{./sections/sec-acknowledgements.tex}
\newpage
\appendix
\input{./sections/sec-coefficients.tex}

\nocite{*}
\newpage
\printbibliography

\end{document}

%% file: preamble.tex
\usepackage{amsthm}
\theoremstyle{definition}

\newtheorem{remark}{Remark}
\newtheorem{example}{Example}
\newtheorem{step}{Step}

\usepackage{amsmath}
\usepackage{physics2}
\usephysicsmodule{ab, ab.legacy}
\usephysicsmodule{qtext.legacy}
\usephysicsmodule{op.legacy}
\usephysicsmodule{nabla.legacy}
\usepackage{derivative, fixdif}

\usepackage[colorlinks,linkcolor=blue,citecolor=blue,urlcolor=blue]{hyperref}
\usepackage{cleveref}
\crefname{step}{step}{step}
\crefname{equation}{Eq.}{Eq.}
\crefname{figure}{Fig.}{Fig.}
\crefname{table}{Tab.}{Tab.}
\crefname{example}{Example.}{Example.}

\numberwithin{equation}{section}

\usepackage{indentfirst}

\usepackage{amsfonts}
\usepackage{mathrsfs}

\usepackage[exponent-mode = scientific]{siunitx}

\def\half{\frac{1}{2}}
\def\dx{\Delta x}
\def\dy{\Delta y}
\def\dt{\Delta t}
\def\tri{\triangle}

\usepackage{bm}
\def\vec#1{\bm{#1}}

\usepackage{graphicx}
\usepackage{caption}
\usepackage{subcaption}
\usepackage{float} 

\usepackage{booktabs} 

\usepackage{tikz}

\usepackage[giveninits=true,
doi=false,
       url=false,
          eprint=false,
         isbn=false,
			backend=biber,
			style=numeric,
			sorting=nyt,
			sortcites=true,
			maxnames=5,
		    minnames=1,
	maxcitenames=2,   
			mincitenames=1,   
]{biblatex}
\DeclareDelimFormat[textcite,parencite]{nametitledelim}{\addcomma\space}
\DeclareDelimFormat[bib]{nametitledelim}{\addcomma\space}
\AtEveryBibitem{
  \clearfield{number}  
  \clearfield{month}  
  \clearfield{day}  
  \clearfield{issue}  
}
\DeclareFieldFormat[article]{title}{#1}
\DeclareFieldFormat[misc]{title}{#1}
\DeclareFieldFormat[inproceedings]{title}{#1}
\DeclareFieldFormat{booktitle}{#1}  
\DeclareFieldFormat[online]{url}{#1}
\DeclareFieldFormat[article]{journaltitle}{\textup{#1},}

\renewbibmacro{in:}{,}

\usepackage{multirow}

%% file: sections/sec-abstract.tex
\section*{Abstract}

In this paper, we propose a derivative-based finite volume Hermite WENO (HWENO) scheme for hyperbolic conservation laws, where both the solution and its first-order derivatives are evolved in time and utilized in spatial reconstructions. The key challenge for solving hyperbolic conservation laws is the possible emergence of discontinuities in the numerical solutions. When facing  discontinuities, the derivatives can become excessively large, which may compromise the robustness of HWENO schemes.  In the first HWENO scheme  \cite{Qiu2004HermiteWENOa}, different sets of stencils were adopted for reconstructing the governing equation and the derivative equation, respectively, aiming to reduce the influence of the derivatives while preserving high-order accuracy. However, this approach not only substantially increases computational cost but also introduces considerable algorithmic complexity. To overcome these limitations, we exclude the information of the target cell's derivatives from spatial reconstructions, while  employing the same reconstructed polynomial during temporal evolution to limit the derivatives. This strategy enhances the robustness of traditional HWENO schemes and allows unified stencils within the derivative-based HWENO framework. Furthermore, the proposed scheme supports arbitrary positive linear weights that sum to one and maintains a compact stencil. Numerical results demonstrate the high-order accuracy, efficiency, high resolution, and robustness of the proposed HWENO scheme.

\noindent \textbf{Keywords:} Hyperbolic conservation laws, Finite volume scheme, HWENO scheme, Derivative-based, High order.

%% file: sections/sec-introduction.tex
\section{Introduction}
In this paper, a fifth-order finite volume Hermite weighted essentially non-oscillatory (Hermite WENO, HWENO) scheme is constructed for both one-dimensional and two-dimensional hyperbolic conservation laws modeled by the following equations:
\begin{equation}
	\left\{
	\begin{aligned}
		 & u_{t} + \div F(u) = 0, \quad          \\
		 & u(\vec{x}, 0) = u_{0}(\vec{x}). \quad \\
	\end{aligned}
	\right.
	\label{eq:hyperbolic-conservation-laws}
\end{equation}
The proposed HWENO scheme introduces derivative equations, but employs a unified set of reconstruction stencils to discretize both the governing equation and the derivative equations, which constitutes the essential distinction from the first HWENO scheme \cite{Qiu2004HermiteWENOa}.
HWENO schemes evolved from the weighted essentially non-oscillatory (WENO) schemes, which have been widely used over the past three decades.
The first WENO scheme was proposed by \textcite{Liu1994WeightedEssentiallya} in 1994,
where all candidate stencils of essentially non-oscillatory (ENO) schemes
\cite{Harten1987UniformlyHighOrdera}
are combined to achieve a third-order accuracy finite volume scheme. In \citeyear{Jiang1996EfficientImplementationa} \cite{Jiang1996EfficientImplementationa}, Jiang and Shu proposed an efficient implementation approach to the WENO scheme in the finite difference framework.
They gave a general
definition of smoothness indicators and nonlinear weights. After that, the WENO schemes boomed with various improvements and extensions
in \cite{Borges2008ImprovedWeighted, Balsara2016EfficientClass, Castro2011HighOrder, Hu1999WeightedEssentially, Levy1999CentralWENO}. A mapping function that maps the nonlinear weights to the linear weights
is designed in \cite{Henrick2005MappedWeighted}
and an improved definition of nonlinear weights is designed in \cite{Castro2011HighOrder}
to achieve optimal order near critical points.
In \cite{Hu1999WeightedEssentially} the WENO scheme is implemented on the triangular meshes.
In \cite{Levy1999CentralWENO} Central WENO schemes are proposed for hyperbolic conservation
laws.
A recent review can be found in \cite{Shu2020EssentiallyNonoscillatory}.
In \cite{Zhu2016NewFifth, Zhu2017NewType}
a novel idea of using artificial linear weights is proposed, which avoids the complicated procedure of computing linear weights
in the previous WENO schemes. It makes the WENO scheme easier to implement
and exhibits better performance in both efficiency and resolution.

The primary difference between HWENO and WENO schemes lies in the fact that the former introduces additional variables, such as derivatives, during spatial reconstruction. Consequently, it can achieve the same or even higher order of accuracy compared to the WENO scheme while using a smaller reconstruction stencil. However, since the solutions of hyperbolic conservation laws may contain discontinuities, the values of these derivatives can become extremely large. As a result, nearly all existing HWENO schemes incorporate some form of treatment for the derivatives. For example, the first HWENO scheme proposed by Qiu and Shu in \cite{Qiu2004HermiteWENOa} employed different sets of reconstruction stencils for the governing equation and the derivative equation in the spatial discretization. This idea was later extended to two-dimensional problems \cite{Qiu2005HermiteWENOb}. However, its robustness in two dimensions is relatively weak.
In subsequent work \cite{Zhu2008ClassFourth}, the robustness was improved by modifying the reconstruction stencil for the derivative equation. Nevertheless, selecting an appropriate reconstruction stencil remains a highly experience-dependent issue. This idea has since been extended to various frameworks, including the finite difference HWENO scheme \cite{Liu2015FiniteDifference}, the moment-based central HWENO scheme \cite{Tao2016HighorderCentral}, the positivity-preserving HWENO scheme \cite{Cai2016PositivityPreservingHigh}, the compact HWENO scheme \cite{Ma2018HWENOSchemes}, the seventh-order HWENO scheme \cite{Zahran2016SeventhOrder}, as well as different types of equations, such as the Hamilton-Jacobi equation \cite{Qiu2005HermiteWENOc}, the Vlasov equation \cite{Yang2014ConservativeNonconservative}, and the KdV equations \cite{Luo2016HybridLDGHWENO}. However, this approach of selecting different reconstruction stencils not only increases the algorithmic complexity and computational cost, but also leaves the first-order derivatives or first-order moments near discontinuities still being used without any restriction, which consequently compromises the robustness of the HWENO scheme.

To improve the robustness of HWENO schemes, \textcite{Zhao2020HybridHermite} drew inspiration from the limiter techniques in the discontinuous Galerkin (DG) methods \cite{Cockburn1989TVBRungeKutta} and proposed a strategy of modifying the first-order moments before performing the spatial reconstruction. This approach effectively enhances the robustness of the scheme, as demonstrated by achieving Courant-Friedrichs-Lewy (CFL) numbers comparable to those of WENO schemes. Subsequently, this work was further extended, by utilizing the WENO reconstruction \cite{Zhu2016NewFifth}, to scenarios allowing arbitrary positive linear weights that sum to one \cite{Zhao2020HermiteWENO}.
In addition, the idea of limiting the first-order moments has also been extended to the finite difference HWENO scheme \cite{Zhao2020ModifiedFiftha} and the Hermite targeted ENO scheme \cite{Wibisono2021FifthOrderHermite}. However, this approach does not resolve the issue of using two different sets of reconstruction stencils in the traditional HWENO schemes \cite{Qiu2004HermiteWENOa,Qiu2005HermiteWENOb}, as it still employs distinct stencils for the spatial reconstruction and the first-order moment limiting procedure in \cite{Zhao2020HybridHermite}. To achieve a unified set of reconstruction stencils, \textcite{Li2021MultiresolutionHWENO} designed multi-resolution finite difference and finite volume HWENO reconstructions based on function values and derivative values. Their approach primarily controls the derivatives through the design of hierarchical stencils and the assignment of corresponding weights. However, this method cannot fully control the derivative values near strong shock waves. For instance, when simulating the forward step problem, their sixth-order finite volume HWENO scheme performed worse than the fourth-order finite difference HWENO scheme. This observation was noted in the comments of their paper, where the underlying reason was stated as unclear. Recently, \textcite{Fan2025MomentBasedHermite} integrated the modification of the first-order moments with the spatial reconstruction into a unified framework, where the modified first-order moments are only used in the temporal discretization. This led to the development of a moment-based HWENO scheme with a unified set of stencils. Throughout the entire process, only a single set of reconstruction stencils is employed. While preserving the robustness of the original moment-limiting HWENO scheme \cite{Zhao2020HermiteWENO}, this approach further simplifies the algorithm and reduces computational cost. This work has also been further extended to the triangular mesh \cite{Zhao2025HighorderMomentbased} and finite difference framework \cite{Xie2025EfficientFifthorder}.

Overall, within the finite volume framework, among the derivative-based HWENO schemes originating from the first HWENO formulation \cite{Qiu2004HermiteWENOa} that reconstructs both function values and derivatives, only the HWENO scheme by \textcite{Li2021MultiresolutionHWENO} has so far achieved a unified set of reconstruction stencils. However, it did not deliver the expected good performance when faced with strong shock wave problems. Therefore, we aim to develop a robust derivative-based HWENO scheme that uses only a single set of reconstruction stencils, thereby addressing the issues inherent in the first HWENO scheme \cite{Qiu2004HermiteWENOa}, such as insufficient robustness and the use of two different sets of reconstruction stencils. We draw inspiration from the moment-based HWENO scheme \cite{Zhao2020ModifiedFiftha}, where the information of the target cell is not used in the spatial reconstruction. Specifically, we reconstruct a nonlinear fifth-degree polynomial. This polynomial is then utilized both for obtaining point values in the spatial discretization and for modifying the derivatives in the temporal discretization. As a result, only a single set of reconstruction stencils is employed throughout the entire process, making full use of the information provided by the spatial reconstruction and improving computational efficiency. Furthermore, in the spatial reconstruction, we can also use arbitrary positive linear weights that sum to one just as in the WENO scheme \cite{Zhu2017NewType}. Compared to the WENO scheme \cite{Zhu2017NewType}, the proposed HWENO scheme not only features a more compact reconstruction stencil but also achieves higher computational efficiency. At the same time, the implementation of a unified stencil significantly improves the computational efficiency of traditional HWENO schemes, requiring less computation time and yielding smaller numerical errors.

This paper is organized as follows. In \cref{sec:description-of-hweno-scheme}, we present the details of the implementation of the derivative-based HWENO scheme in both one- and two-dimensional cases. In \cref{sec:numerical-tests}, several benchmark numerical examples, including one-dimensional and two-dimensional accuracy tests, are conducted to verify the fifth-order accuracy, efficiency, and high resolution of the constructed scheme. Finally, some concluding remarks are given in \cref{sec:conclusions}.

%% file: sections/sec-description.tex
\section{Description of HWENO Scheme}
\label{sec:description-of-hweno-scheme}
In this section, we present the construction methodologies of the derivative-based fifth-order finite volume HWENO schemes
for both one-dimensional and two-dimensional hyperbolic conservation laws on uniform meshes.
\input{./sections/subsec-one-dimensional-case.tex}
\input{./sections/subsec-uniform-case.tex}

%% file: sections/subsec-one-dimensional-case.tex
\subsection{HWENO Scheme in One-Dimensional Case}
Consider one-dimensional scalar hyperbolic conservation laws,
\begin{equation}
	\left\{
	\begin{aligned}
		 & u_{t} + f(u)_{x} = 0, \quad \\
		 & u(x, 0) = u_{0}(x) . \quad  \\
	\end{aligned}
	\right.
	\label{eq:1D-scalar-conservation-laws}
\end{equation}
For simplicity, the computational domain is divided into cells $ I_{i} = [x_{i-\half}, x_{i+\half}] $,
with uniform cell size $ \dx  = x_{i+\half} - x_{i-\half}  $, and the cell centers are $ x_{i} = \frac{1}{2}(x_{i-\half} + x_{i+\half}) $.
Taking the derivative of \cref{eq:1D-scalar-conservation-laws} with respect to $ x $, we have
\begin{equation*}
	\left\{
	\begin{aligned}
		 & u_{t} + f(u)_{x} = 0, \quad    \\
		 & v_{t} + h(u, v)_{x} = 0, \quad \\
	\end{aligned}
	\right.
\end{equation*}
where $ v = u_{x} $ and $ h(u, v) = f'(u)\cdot v $.
Integrating over the cell $ I_{i} $, we have
\begin{equation*}
	\left\{
	\begin{aligned}
		 & \odv*{\bar{u}_{i}(t)}{t} = -\frac{1}{\dx }[f(u(x_{i-\half},t)) - f(u(x_{i+\half},t))]  , \quad                                      \\
		 & \odv*{\bar{v}_{i}(t)}{t} = -\frac{1}{\dx }[h(u(x_{i-\half},t), v(x_{i-\half},t)) - h(u(x_{i+\half},t), v(x_{i+\half},t))]  , \quad \\
	\end{aligned}
	\right.
\end{equation*}
where $ \bar{u}_{i}(t) := \frac{1}{\dx }\int_{I_{i}}^{} u(x, t) \d{x}$ and
$ \bar{v}_{i}(t) := \frac{1}{\dx }\int_{I_{i}}^{} v(x, t) \d{x}$
are the cell averages of $ u(x, t) $ and $ v(x, t) $ over the cell $ I_{i} $, respectively.

By approximating the value of the flux functions $ f(u(x_{i+\half},t)) $ and $ h(u(x_{i+\half},t), v(x_{i+\half},t)) $
by numerical fluxes $ \hat{f}_{i+\half} = \hat{f}(u_{i+\half}^{-}, u_{i+\half}^{+}) $ and $ \hat{h}_{i+\half} = \hat{h}(u_{i+\half}^{-}, u_{i+\half}^{+}, v_{i+\half}^{-}, v_{i+\half}^{+}) $, respectively,
we obtain a conservative semi-discrete scheme
\begin{equation}
	\left\{
	\begin{aligned}
		 & \odv*{\bar{u}_{i}}{t} = -\frac{1}{\dx }(\hat{f}_{i+\half} - \hat{f}_{i-\half}) =: L_{i}^{u}, \quad \\
		 & \odv*{\bar{v}_{i}}{t} = -\frac{1}{\dx }(\hat{h}_{i+\half} - \hat{h}_{i-\half}) =: L_{i}^{v}, \quad \\
	\end{aligned}
	\right.
	\label{eq:semi-discrete-scheme}
\end{equation}
where $ \bar{u}_{i} $ and $ \bar{v}_{i} $ are the numerical approximations of $ \bar{u}_{i}(t) $ and $ \bar{v}_{i}(t) $, respectively.
The fluxes $ \hat{f}_{i+\half} $ and $ \hat{h}_{i+\half} $ are approximated by the global Lax-Friedrichs fluxes
\begin{equation}
	\left\{
	\begin{aligned}
		\hat{f}_{i+\half} & = \frac{1}{2}[f(u_{i+\half}^{-}) + f(u_{i+\half}^{+})] - \frac{\alpha}{2} (u_{i+\half}^{+} - u_{i+\half}^{-}),                                   \\
		\hat{h}_{i+\half} & = \frac{1}{2}[h(u_{i+\half}^{-}, v_{i+\half}^{-}) + h(u_{i+\half}^{+}, v_{i+\half}^{+})] - \frac{\alpha}{2} (v_{i+\half}^{+} - v_{i+\half}^{-}),
	\end{aligned}
	\right.
\end{equation}
where $ \alpha = \max_{u} \abs{f'(u)} $.
$ u_{i+\half}^{-} $ and $ u_{i+\half}^{+} $ denote the left and right limits of numerical solution $ u(x) $ at the cell interface $ x_{i+\half} $, respectively.
$ v_{i+\half}^{-} $ and $ v_{i+\half}^{+} $ denote the left and right limits of the derivative of the numerical solution $ u(x) $ at the cell interface $ x_{i+\half} $, respectively.

Now we describe the reconstruction procedure of $ u_{i+\half}^{-} $, $ u_{i-\half}^{+} $, $ v_{i+\half}^{-} $ and $ v_{i-\half}^{+} $,
and present the time discretization method.
\begin{step}
	Reconstruct a quartic polynomial $ p_{0}(x) $, and two linear ones $ p_{1}(x) $ and $ p_{2}(x) $.
\end{step}
Firstly, a quartic polynomial $ p_{0}(x) = \sum_{\ell = 0}^{4} c_{0,\ell} (\frac{x - x_{i}}{\dx })^{\ell}  $
is reconstructed on a big stencil $ S_{0} = \{I_{i-1}, I_{i}, I_{i+1} \} $,
satisfying
\begin{equation}
	\left\{
	\begin{aligned}
		 & \frac{1}{\dx }\int_{I_{k}}^{} p_{0}(x) \d{x} = \bar{u}_{k},\quad k = i-1, i, i+1 \\
		 & \frac{1}{\dx }\int_{I_{k}}^{} p_{0}'(x) \d{x} = \bar{v}_{k},\quad k = i-1,  i+1. \\
	\end{aligned}
	\right.
\end{equation}
Secondly, two linear polynomials $ p_{1}(x) = \sum_{\ell=0}^{1} c_{1,\ell}(\frac{x - x_{i}}{\dx })^{\ell}$
and $ p_{2}(x) = \sum_{\ell=0}^{1} c_{2,\ell}(\frac{x - x_{i}}{\dx })^{\ell}   $
are reconstructed on stencil $ S_{1} = \{I_{i-1}, I_{i}\} $
and $ S_{2} = \{I_{i}, I_{i+1}\} $, respectively, satisfying
\begin{equation}
	\left\{
	\begin{aligned}
		 & \frac{1}{\dx }\int_{I_{k}}^{} p_{1}(x) \d{x} = \bar{u}_{k},\quad k = i-1, i, \\
		 & \frac{1}{\dx }\int_{I_{k}}^{} p_{2}(x) \d{x} = \bar{u}_{k},\quad k = i, i+1. \\
	\end{aligned}
	\right.
\end{equation}
The explicit expressions of the coefficients $ c_{m,\ell} $ are given in Appendix~\ref{sec:coefficients}.
Finally, we rewrite $ p_{0}(x) $ as
\begin{equation}
	p_{0}(x) = \gamma_{0} (\frac{1}{\gamma _{0}}p_{0}(x) - \sum_{m = 1}^{2} \frac{\gamma_{m}}{\gamma _{0}}p_{m}(x))
	+ \sum_{m = 1}^{2} \gamma_{m} p_{m}(x),
	\label{eq:rewrite}
\end{equation}
where $ \gamma_{m} $ is positive for $ m = 0, 1, 2 $ and satisfy $ \sum_{m = 0}^{2} \gamma_{m} = 1 $.
\begin{step}
	Compute the smoothness indicators $ \{\beta _{m}\}_{m = 0}^{2} $
	to measure the smoothness level for each polynomial in cell $ I_{i} $.
\end{step}
The smoothness indicator in the classical WENO scheme \cite{Jiang1996EfficientImplementationa} is defined as
\begin{equation}
	\beta^{\text{typical}}  = \sum_{\ell = 1}^{\deg p(x)}\frac{1}{\dx } \int_{I_{i}}^{}\dx ^{2\ell} [\odv*[\ell]{p(x)}{x}]^{2}  \d{x},
	\label{eq:smoothness-indicator-typical}
\end{equation}
To reduce computational cost, we employ a new definition of smoothness indicator \cite{Zhao2025HighorderMomentbased} for any polynomial
$ p(x) = \sum_{\ell = 0}^{\deg p(x)}c_{\ell} (\frac{x - x_{i}}{\dx })^{\ell}  $
\begin{equation}
	\beta = \sum_{\ell = 1}^{\deg p(x)} [\odv*[\ell]{p(x_{i})}{x}]^{2}= \sum_{\ell = 1}^{\deg p(x)} (\ell!)^{2} c_{\ell},
\end{equation}
where $ c_{m, \ell} $ are the coefficients of the polynomials $ \{p_{m}(x)\} $ mentioned above.
It is obtained by approximating those integrals in \cref{eq:smoothness-indicator-typical}
by midpoint quadrature formula.

	In \cite{Qiu2004HermiteWENOa}, they mentioned that the summation in the smoothness indicator \cref{eq:smoothness-indicator-typical} for reconstruction of derivatives should start from $ \ell = 2 $,
	however, our schemes do not need special treatment of smoothness indicators for derivative reconstruction.

\begin{step}
	Reconstruct the solution $ u_{i}(x) $ and modify the derivative $ \hat{v}_{i} $ used in time discretization.
	\label{step:HWENO-reconstruction}
\end{step}
In this step, we employ the nonlinear fifth-order HWENO reconstructions for $ u_{i}(x) $,
and modify those average derivatives $ \hat{v}_{i} $ simultaneously.
The reconstruction procedure maintains high-order accuracy in smooth areas
while ensuring essentially non-oscillatory behavior near discontinuities.
As in the WENO and HWENO schemes \cite{Fan2025MomentBasedHermite, Zhu2016NewFifth, Zhu2017NewType, Zhao2025HighorderMomentbased, Zhu2018NewFinite},
a parameter $ \tau $ formulated as
\begin{equation}
	\tau = \frac{\sum_{m = 1}^{2} \abs{\beta _{0} - \beta _{m}}}{2},
\end{equation}
is used to measure the overall difference between $ \{\beta _{m}\}_{m=1}^{2} $ and $ \beta _{0} $,
then the  nonlinear weights are computed by
\begin{equation}
	\omega_{m} = \frac{\tilde{\omega}_{m}}{\sum_{m = 0}^{2} \tilde{\omega}_{m} }, \;
	\text{with } \tilde{\omega}_{m} = \gamma _{m}(1 + \frac{\tau^{2}}{\beta_{m} + \varepsilon}), \;
	\text{for } m = 0, 1, 2,
	\label{eq:nonlinear-weights}
\end{equation}
where $ \varepsilon = 10^{-8} $ is a small positive number to avoid division by zero.

Finally, a nonlinear HWENO reconstruction of the polynomial $ u_{i}(x) $ for $ u(x) $
within target cell $ I_{i} $ is obtained by replacing some of the linear weights in
\cref{eq:rewrite} by nonlinear weights defined in \cref{eq:nonlinear-weights}
\begin{equation}
	u_{i}(x) = \omega_{0} (\frac{1}{\gamma _{0}}p_{0}(x) - \sum_{m = 1}^{2} \frac{\gamma_{m}}{\gamma _{0}}p_{m}(x))
	+ \sum_{m = 1}^{2} \omega_{m} p_{m}(x).
	\label{eq:hweno-reconstruction-1D}
\end{equation}
Concurrently, the modified average derivative $ \hat{v}_{i} $
is determined by averaging the derivative function of $ u_{i}(x) $ over the target cell $ I_{i} $
\begin{equation}
	\hat{v}_{i} = \frac{1}{\dx }\int_{I_{i}}^{} u'_{i}(x) \d{x} = \frac{u_{i}(x_{i+\half}) - u_{i}(x_{i-\half})}{\dx },
\end{equation}
which is only used in the time discretization, inspired by \cite{Fan2025MomentBasedHermite, Zhao2025HighorderMomentbased}.
The required values at interfaces are evaluated by
\begin{equation}
	u_{i-\half}^{+} = u_{i}(x_{i-\half}),\;u_{i+\half}^{-} = u_{i}(x_{i+\half}),\;
	v_{i-\half}^{+} = u'_{i}(x_{i-\half}),\;v_{i+\half}^{-} = u'_{i}(x_{i+\half}).
\end{equation}
\begin{step}
	Time discretization method for the semi-discrete scheme \cref{eq:semi-discrete-scheme}.
\end{step}
To ensure the stability of the scheme, we apply the third-order strong stability preserving (SSP) Runge-Kutta time discretization method \cite{Gottlieb2001StrongStabilityPreserving}.
Due to the inherent presence of discontinuities in solutions to hyperbolic conservation laws,
the magnitude of the derivatives can escalate significantly in non-smooth regions. Therefore, we modify the average derivatives $ \hat{v}_{i} $ in each stage of the SSP Runge-Kutta method.

\begin{equation}
	\left\{
	\begin{aligned}
		\vec{u}^{(1)}_{i}   & = \hat{\vec{u}}^{n}_{i} + \dt L_{i}(\vec{u}^{n}),                                                   \\
		\vec{u}^{(2)}_{i}   & = \frac{3}{4}\hat{\vec{u}}^{n}_{i} + \frac{1}{4}( \hat{\vec{u}}^{(1)}_{i} + \dt L_{i}(\vec{u}^{(1)}) ), \\
		\vec{u}^{n+1}_{i} & = \frac{1}{3}\hat{\vec{u}}^{n}_{i} + \frac{2}{3}( \hat{\vec{u}}^{(2)}_{i} + \dt L_{i}(\vec{u}^{(2)}) ), \\
	\end{aligned}
	\right.
\end{equation}
where $ \vec{u}^{\star}_{i} := [\bar{u}_{i}^{\star}, \bar{v}_{i}^{\star}] $,
$ \hat{\vec{u}}^{\star}_{i} := [\bar{u}_{i}^{\star}, \hat{v}_{i}^{\star}] $,
and $ L_{i} := [L_{i}^{u}, L_{i}^{v}] $.
These modified average partial derivatives are defined as
$ \hat{v}_{i}^{\star} = \frac{1}{\dx }\int_{I_{i}}^{} \pdv*{}{x} u(x,t^{\star}) \d{x} $
where
$ \star $ denotes $ n, (1) $, and $ (2) $ corresponding to each time stage.

%% file: sections/subsec-uniform-case.tex
\subsection{HWENO Scheme in the Two-Dimensional Case}
Consider two-dimensional hyperbolic conservation laws, as defined by
\begin{equation}
	\left\{
	\begin{aligned}
		 & u_{t} + \div F(u) = 0 , \quad  \\
		 & u(x, y,0) = u_{0}(x, y), \quad \\
	\end{aligned}
	\right.
	\label{eq:2D-scalar-conservation-laws}
\end{equation}
where $ F = (f, g) $, $ \div F(u) = f(u)_{x} + g(u)_{y}  $.
Taking the partial derivative of \cref{eq:2D-scalar-conservation-laws} with respect to $ x $ and $ y $, respectively, we have

\begin{equation}
	\left\{
	\begin{aligned}
		 & {u}_{t} + f(u)_{x} + g(u)_{y} = 0, \quad             \\
		 & {v}_{t} + f^{v}(u,v)_{x} + g^{v}(u,v)_{y} = 0, \quad \\
		 & {w}_{t} + f^{w}(u,w)_{x} + g^{w}(u,w)_{y} = 0. \quad \\
	\end{aligned}
	\right.
	\label{eq:2D-uvw-cartesian}
\end{equation}
where $ v = u_{x} $, $ w = u_{y} $, $ f^{v}(u,v) = f'(u) \cdot v $, $ g^{v}(u,v) = g'(u) \cdot v $, $ f^{w}(u,w) = f'(u) \cdot w $, and $ g^{w}(u,w) = g'(u) \cdot w $
.
For simplicity, the computational domain $ \Omega = [a, b] \times [c, d] $ is divided uniformly into $ N_x \times N_y $ rectangular cells $ I_{i,j} = [x_{i-\half}, x_{i+\half}] \times [y_{j-\half}, y_{j+\half}] $ for $ i = 1, 2, \ldots, N_x $, $ j = 1, 2, \ldots, N_y $. The mesh sizes are $ \dx = x_{i+\half} - x_{i-\half} $ and $ \dy = y_{j+\half} - y_{j-\half} $, where $ (x_{i}, y_{j}) $ is the center of cell $ I_{i,j} $, $ x_{i} = \frac{1}{2}(x_{i-\half}+x_{i+\half}) $, $ y_{j} = \frac{1}{2}(y_{j-\half}+y_{j+\half}) $.

Integrating \cref{eq:2D-uvw-cartesian} over the target cell $ I_{i,j} $, we have

\begin{equation*}
	\left\{
	\begin{aligned}
		\odv*{\bar{u}_{i,j}}{t} & = -\frac{1}{\dx\dy}(
		\int_{y_{j-\half}}^{y_{j+\half}} [f(u)]_{x_{i-\half}}^{x_{i+\half}} \d{y}
		+ \int_{x_{i-\half}}^{x_{i+\half}} [g(u)]_{y_{j-\half}}^{y_{j+\half}} \d{x}),
		\\
		\odv*{\bar{v}_{i,j}}{t} & = -\frac{1}{\dy}(
		\int_{y_{j-\half}}^{y_{j+\half}} [f^{v}(u,v)]_{x_{i-\half}}^{x_{i+\half}} \d{y}
		+ \int_{x_{i-\half}}^{x_{i+\half}} [g^{v}(u,v)]_{y_{j-\half}}^{y_{j+\half}} \d{x}),
		\\
		\odv*{\bar{w}_{i,j}}{t} & = -\frac{1}{\dx}(
		\int_{y_{j-\half}}^{y_{j+\half}} [f^{w}(u,w)]_{x_{i-\half}}^{x_{i+\half}} \d{y}
		+ \int_{x_{i-\half}}^{x_{i+\half}} [g^{w}(u,w)]_{y_{j-\half}}^{y_{j+\half}} \d{x}),
	\end{aligned}
	\right.
\end{equation*}
where  $  \bar{u}_{i,j} := \frac{1}{\dx\dy} \int_{I_{i,j}}^{} u \d{x}\d{y} $, $ \bar{v}_{i,j} := \frac{1}{\dy} \int_{I_{i,j}}^{} u_{x} \d{x}\d{y} $ and $ \bar{w}_{i,j}:= \frac{1}{\dx} \int_{I_{i,j}}^{} u_{y} \d{x}\d{y} $ denote the cell average of $ u $ and the scaled cell averages associated with $ u_{x} $ and $ u_{y} $ over the cell $ I_{i,j} $, respectively.
Moreover, $ [h(x)]_{a}^{b} := h(b)-h(a) $ for any function $ h $.
To achieve high-order accuracy, we approximate those integrals of fluxes by the Gaussian quadrature method of sufficient accuracy. To achieve fifth-order accuracy,
one can employ a three-point Gaussian quadrature rule.
Let $ \{x_{i}^{\mu}\}_{\mu = 1}^{Q} $ and
$ \{y_{j}^{\mu}\}_{\mu = 1}^{Q} $
denote the roots of Legendre polynomial in
$ [x_{i-\half}, x_{i+\half}] $ and $ [y_{j-\half}, y_{j+\half}] $, respectively,
and $ \{\sigma_\mu\}_{\mu=1}^{Q} $ denote the corresponding weights,
then we obtain the semi-discrete scheme

\begin{equation}
	\left\{
	\begin{aligned}
		\odv*{\bar{u}_{i,j}}{t} & = -\frac{1}{\dx\dy}(\hat{f}_{i+\half,j} - \hat{f}_{i-\half,j} + \hat{g}_{i,j+\half} - \hat{g}_{i,j-\half}) =:L_{i,j}^{u},
		\\
		\odv*{\bar{v}_{i,j}}{t} & = -\frac{1}{\dy}(\hat{f}^{v}_{i+\half,j} - \hat{f}^{v}_{i-\half,j} + \hat{g}^{v}_{i,j+\half} - \hat{g}^{v}_{i,j-\half})=:L_{i,j}^{v},
		\\
		\odv*{\bar{w}_{i,j}}{t} & = -\frac{1}{\dx}(\hat{f}^{w}_{i+\half,j} - \hat{f}^{w}_{i-\half,j} + \hat{g}^{w}_{i,j+\half} - \hat{g}^{w}_{i,j-\half})=:L_{i,j}^{w},
	\end{aligned}
	\right.
	\label{eq:2D-semi-discrete-scheme-uniform}
\end{equation}
where $ \hat{f}_{i+\half,j} $, $ \hat{g}_{i,j+\half} $, $ \hat{f}^{v}_{i+\half,j} $, $ \hat{g}^{v}_{i,j+\half} $, $ \hat{f}^{w}_{i+\half,j} $, and $ \hat{g}^{w}_{i,j+\half} $
denote the approximations to\\
$ \int_{y_{j-\half}}^{y_{j+\half}} f(u(x_{i+\half},y,t)) \d{y} $, $ \int_{x_{i-\half}}^{x_{i+\half}} g(u(x,y_{j+\half},t)) \d{x} $,
$ \int_{y_{j-\half}}^{y_{j+\half}} f'(u(x_{i+\half},y,t))\cdot v(x_{i+\half},y,t) \d{y} $, $ \int_{x_{i-\half}}^{x_{i+\half}} g'(u(x,y_{j+\half},t))\cdot v(x,y_{j+\half},t) \d{x} $,
$ \int_{y_{j-\half}}^{y_{j+\half}} f'(u(x_{i+\half},y,t))\cdot w(x_{i+\half},y,t) \d{y} $,\\ and $ \int_{x_{i-\half}}^{x_{i+\half}} g'(u(x,y_{j+\half},t))\cdot w(x,y_{j+\half},t) \d{x} $, respectively.
And the numerical fluxes are defined as
\begin{equation}
	\left\{
	\begin{aligned}
		 & \hat{f}_{i+\half,j} = \dy\sum_{\mu=1}^{Q} \sigma_{\mu} \hat{f}(u_{i+\half,j}^{-,\mu},u_{i+\half,j}^{+,\mu}),                                                     \\
		 & \hat{g}_{i,j+\half} = \dx\sum_{\mu=1}^{Q} \sigma_{\mu} \hat{g}(u_{i,j+\half}^{\mu,-},u_{i,j+\half}^{\mu,+}),                                                     \\
		 & \hat{f}^{v}_{i+\half,j} = \dy\sum_{\mu=1}^{Q} \sigma_{\mu} \hat{f}^{v}(u_{i+\half,j}^{-,\mu},u_{i+\half,j}^{+,\mu},v_{i+\half,j}^{-,\mu},v_{i+\half,j}^{+,\mu})  \\
		 & \hat{g}^{v}_{i,j+\half} = \dx\sum_{\mu=1}^{Q} \sigma_{\mu} \hat{g}^{v}(u_{i,j+\half}^{\mu,-},u_{i,j+\half}^{\mu,+},v_{i,j+\half}^{\mu,-},v_{i,j+\half}^{\mu,+}), \\
		 & \hat{f}^{w}_{i+\half,j} = \dy\sum_{\mu=1}^{Q} \sigma_{\mu} \hat{f}^{w}(u_{i+\half,j}^{-,\mu},u_{i+\half,j}^{+,\mu},w_{i+\half,j}^{-,\mu},w_{i+\half,j}^{+,\mu})  \\
		 & \hat{g}^{w}_{i,j+\half} = \dx\sum_{\mu=1}^{Q} \sigma_{\mu} \hat{g}^{w}(u_{i,j+\half}^{\mu,-},u_{i,j+\half}^{\mu,+},w_{i,j+\half}^{\mu,-},w_{i,j+\half}^{\mu,+}), \\
	\end{aligned}
	\right.
	\nonumber
\end{equation}
with the global Lax-Friedrichs fluxes
\begin{equation}
	\left\{
	\begin{aligned}
		\hat{f}(u^{-}, u^{+})                 & = \frac{1}{2}[f(u^{-}) + f(u^{+})] - \frac{\alpha_{x}}{2} (u^{+} - u^{-}),                     \\
		\hat{g}(u^{-}, u^{+})                 & = \frac{1}{2}[g(u^{-}) + g(u^{+})] - \frac{\alpha_{y}}{2} (u^{+} - u^{-}),                     \\
		\hat{f}^{v}(u^{-}, u^{+},v^{-},v^{+}) & = \frac{1}{2}[f^{v}(u^{-},v^{-}) + f^{v}(u^{+},v^{+})] - \frac{\alpha_{x}}{2} (v^{+} - v^{-}), \\
		\hat{g}^{v}(u^{-}, u^{+},v^{-},v^{+}) & = \frac{1}{2}[g^{v}(u^{-},v^{-}) + g^{v}(u^{+},v^{+})] - \frac{\alpha_{y}}{2} (v^{+} - v^{-}), \\
		\hat{f}^{w}(u^{-}, u^{+},w^{-},w^{+}) & = \frac{1}{2}[f^{w}(u^{-},w^{-}) + f^{w}(u^{+},w^{+})] - \frac{\alpha_{x}}{2} (w^{+} - w^{-}), \\
		\hat{g}^{w}(u^{-}, u^{+},w^{-},w^{+}) & = \frac{1}{2}[g^{w}(u^{-},w^{-}) + g^{w}(u^{+},w^{+})] - \frac{\alpha_{y}}{2} (w^{+} - w^{-}), \\
	\end{aligned}
	\right.
\end{equation}
where $ \alpha_{x} = \max_{u} \abs{f'(u)} $ and $ \alpha_{y} = \max_{u} \abs{g'(u)} $.
$ u_{i+\half,j}^{\pm,\mu} $ and $ u_{i,j+\half}^{\mu,\pm} $  denote the
reconstructed values of $ u $ at the Gaussian points $ (x_{i+\half}, y_{j}^{\mu}) $ and $ (x_{i}^{\mu}, y_{j+\half}) $, respectively, and so do $ v $ and $ w $.
The superscripts $ - $ and $ + $ indicate the value is approximated inside and outside the target cell $ I_{i,j} $, respectively.

Now we describe the reconstruction procedure of those reconstructed values at interfaces,
and present the time discretization method.

\begin{figure}[htbp]
	\centering
	\begin{subfigure}{0.49\textwidth}
		\centering
		\begin{tikzpicture}[scale=1.6]
			\draw [thick] (0, 0) grid (3, 3);
			\draw (1.5, 1.5) node {$ I_{i,j} $};
			\draw (1.5-1, 1.5) node {$ I_{i-1,j} $};
			\draw (1.5, 1.5-1) node {$ I_{i,j-1} $};
			\draw (1.5+1, 1.5) node {$ I_{i+1,j} $};
			\draw (1.5, 1.5+1) node {$ I_{i,j+1} $};
			\draw (1.5-1, 1.5-1) node {$ I_{i-1,j-1} $};
			\draw (1.5+1, 1.5-1) node {$ I_{i+1,j-1} $};
			\draw (1.5+1, 1.5+1) node {$ I_{i+1,j+1} $};
			\draw (1.5-1, 1.5+1) node {$ I_{i-1,j+1} $};
		\end{tikzpicture}
		\caption{The stencil $ S_{0} $ for the quartic polynomial.}
	\end{subfigure}
	\begin{subfigure}{0.49\textwidth}
		\centering
		\begin{tikzpicture}[scale=1.1]
			\draw [thick] (0, 0) rectangle (1, 1);
			\draw [thick] (0, 1) rectangle (1, 2);
			\draw [thick] (-1, 0) rectangle (0, 1);
			\draw (.5, .5) node {$ I_{i,j} $};
			\draw (.5-1, .5) node {$ I_{i-1,j} $};
			\draw (.5, .5+1) node {$ I_{i,j+1} $};
			\draw (.5-1, .5+1) node {$ S_4 $};

			\draw [thick, shift={(1.5,0)}] (0, 0) rectangle (1, 1);
			\draw [thick, shift={(1.5,0)}] (1, 0) rectangle (2, 1);
			\draw [thick, shift={(1.5,0)}] (0, 1) rectangle (1, 2);
			\draw [shift={(1.5,0)}] (.5, .5) node {$ I_{i,j} $};
			\draw [shift={(1.5,0)}] (.5+1, .5) node {$ I_{i+1,j} $};
			\draw [shift={(1.5,0)}] (.5, .5+1) node {$ I_{i,j+1} $};
			\draw [shift={(1.5,0)}] (.5+1, .5+1) node {$ S_3 $};

			\draw [thick, shift={(0,-1.5)}] (0, 0) rectangle (1, 1);
			\draw [thick, shift={(0,-1.5)}] (-1, 0) rectangle (0, 1);
			\draw [thick, shift={(0,-1.5)}] (0, -1) rectangle (1, 0);
			\draw [shift={(0,-1.5)}] (.5, .5) node {$ I_{i,j} $};
			\draw [shift={(0,-1.5)}] (.5-1, .5) node {$ I_{i-1,j} $};
			\draw [shift={(0,-1.5)}] (.5, .5-1) node {$ I_{i,j-1} $};
			\draw [shift={(0,-1.5)}] (.5-1, .5-1) node {$ S_{1} $};

			\draw [thick, shift={(1.5,-1.5)}] (0, 0) rectangle (1, 1);
			\draw [thick, shift={(1.5,-1.5)}] (0, 0-1) rectangle (1, 1-1);
			\draw [thick, shift={(1.5,-1.5)}] (0+1, 0) rectangle (1+1, 1);
			\draw [shift={(1.5,-1.5)}] (.5, .5) node {$ I_{i,j} $};
			\draw [shift={(1.5,-1.5)}] (.5+1, .5) node {$ I_{i+1,j} $};
			\draw [shift={(1.5,-1.5)}] (.5, .5-1) node {$ I_{i,j-1} $};
			\draw [shift={(1.5,-1.5)}] (.5+1, .5-1) node {$ S_{2} $};
		\end{tikzpicture}
		\caption{The stencils for the linear polynomials.}
	\end{subfigure}
	\caption{The stencils for two-dimensional HWENO reconstruction.}
	\label{fig:stencils}
\end{figure}
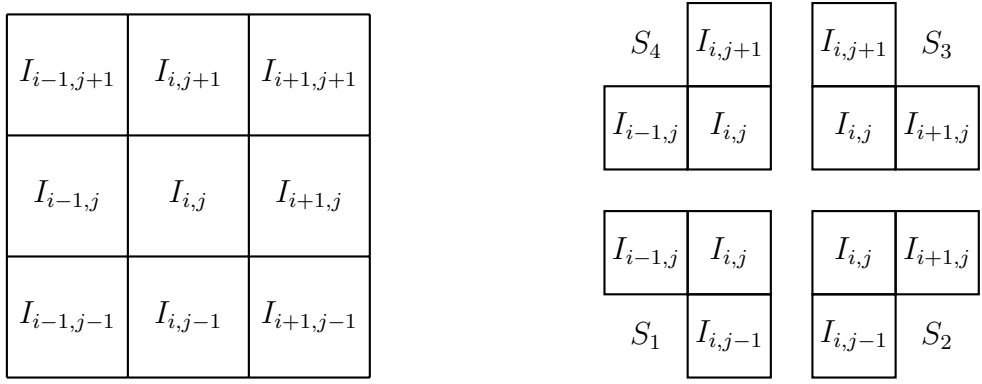

\setcounter{step}{0}
\begin{step}
	Reconstruct a quartic bivariate polynomial $ p_{0}(x, y) $ and four linear ones $ \{p_{m}(x, y)\}_{m = 1}^{4} $.
\end{step}
Let $ \phi $ denote the natural basis functions over the target cell $ I_{i,j} $, defined as
\begin{equation}
	\phi = \cup_{n=0}^{+\infty} \{\xi^{k}\eta^{n-k}:k = 0,1,\ldots,n\},\;
	\text{with } \xi = \frac{x - x_{i}}{\dx},\; \eta = \frac{y - y_{i}}{\dy},
\end{equation}
thus we have $ p_{0}(x, y) = \sum_{\ell = 1}^{15} c_{0,\ell}\phi_{\ell}(x, y)  $ and $ p_{m}(x, y) = \sum_{\ell = 1}^{3} c_{m,\ell}\phi_{\ell}(x,y) $,
for $ m = 1, 2, 3, 4 $.
The quartic polynomial $ p_{0}(x,y) $ is obtained by HWENO reconstruction on a large stencil $ S_{0} = \{I_{i+\alpha,j+\beta}: \alpha = -1,0,1,\, \beta = -1,0,1.\} $,
and four linear ones $ p_{m}(x, y) $ are reconstructed
on four small stencils $ S_{1} = \{I_{i,j},I_{i-1,j},I_{i,j-1}\} $,
$ S_{2} = \{I_{i,j},I_{i+1,j},I_{i,j-1}\} $,
$ S_{3} = \{I_{i,j},I_{i+1,j},I_{i,j+1}\} $,
$ S_{4} = \{I_{i,j},I_{i-1,j},I_{i,j+1}\} $.
\cref{fig:stencils} illustrates the selection of these stencils.

The reconstruction of the quartic polynomial $ p_{0}(x,y) $ requires
\begin{equation}
	\left\{
	\begin{aligned}
		 & \frac{1}{\dx\dy} \int_{I}^{} p_0(x,y) \d{\Omega} = \bar{u}_{I}, \quad \text{for } I \in S_0,                                      \\
		 & \frac{1}{\dy} \int_{I}^{} \pdv*{p_0(x, y)}{x}  \d{\Omega} \approx \bar{v}_{I}  , \quad \text{for } I \in S_0\setminus\{I_{i,j}\}, \\
		 & \frac{1}{\dx} \int_{I}^{} \pdv*{p_0(x, y)}{y}  \d{\Omega} \approx \bar{w}_{I}  , \quad \text{for } I \in S_0\setminus\{I_{i,j}\},
	\end{aligned}
	\right.
\end{equation}
among which the equality constraints are precisely satisfied, while the approximate equality constraints are satisfied in the least squares sense.
And for those linear polynomials $ p_{m}(x, y) $, they are uniquely determined by
\begin{equation}
	\frac{1}{\dx\dy} \int_{I}^{} p_m(x,y) \d{\Omega} = \bar{u}_{I}, \quad \text{for } I \in S_m,\, m = 1,2,3,4.
\end{equation}
The details of the coefficients of these polynomials are given in Appendix~\ref{sec:coefficients}.
Then, we reformulate the polynomial $ p_{0}(x, y) $ as
\begin{equation}
	p_{0}(x,y) = \gamma_{0} (\frac{1}{\gamma _{0}}p_{0}(x,y) - \sum_{m = 1}^{4} \frac{\gamma_{m}}{\gamma _{0}}p_{m}(x,y))
	+ \sum_{m = 1}^{4} \gamma_{m} p_{m}(x),
	\label{eq:rewrite-2D}
\end{equation}
To ensure the stability of the resulting HWENO reconstruction,
the set of linear weights $ \{\gamma_{m}\}_{m = 0}^{4} $ must be positive,
with the requirement $ \sum_{m = 0}^{4} \gamma_{m} = 1  $ again.

\begin{step}
	Compute the smoothness indicators $ \{\beta_{m}\}_{m = 0}^{4} $ for $ \{p_{m}(x, y)\}_{m = 0}^{4} $, which measure the level of smoothness for each polynomial.
\end{step}
The smoothness indicator for a polynomial $ p(x, y) $ is typically defined as
\begin{equation}
	\beta^{\text{typical}}  =
	\sum_{\ell = 1}^{\deg p(x,y)}\frac{1}{\abs{I_{i,j}}}
	\sum_{m+n = \ell}^{}
	\int_{\tri_{i}}^{}  [\dx^{m}\dy^{n}\pdv*[m,n]{p(x, y)}{x,y}]^{2}  \d{\Omega}.
	\label{eq:smoothness-indicator-typical-2D}
\end{equation}
To reduce computational cost, we employ a new
definition \cite{Zhao2025HighorderMomentbased} again
\begin{equation}
	\beta =
	\sum_{\ell = 1}^{\deg p(x,y)}\sum_{m+n = \ell}^{} [\pdv*[m,n]{p(x_{i}, y_{i})}{x,y}]^{2}  =
	\sum_{\ell = 1}^{\deg p(x,y)}\sum_{m+n = \ell}^{} (m!n!)^{2} c_{\star}^{2},
\end{equation}
where $ c_{\star} $ is the corresponding coefficient of the term $ \xi^{m}\eta^{n} $ in the polynomial $ p(x, y) $.

\begin{step}
	Reconstruct the solution $ u_{i}(x, y) $ and obtain the modified partial derivatives $ \hat{v}_{i} $ and $ \hat{w}_{i} $ by HWENO reconstruction.
\end{step}
In this step, we employ the nonlinear high-order HWENO reconstructions for $ u_{i}(x, y) $, $ \hat{v}_{i} $ and $ \hat{w}_{i} $, respectively.
These reconstructions maintain high-order accuracy in smooth areas while essentially ensuring non-oscillatory behavior near discontinuities.
As in the WENO and HWENO schemes, a parameter $ \tau $ is used to measure
the overall difference between the set of $ \{\beta_{m}\}_{m=1}^{4} $ and $ \beta_{0} $, formulated as
\begin{equation}
	\tau = \frac{1}{4} \sum_{m = 1}^{4} \abs{\beta _{0} - \beta _{m}},
\end{equation}
then, the nonlinear weights are computed by
\begin{equation}
	\omega_{m} = \frac{\tilde{\omega}_{m}}{\sum_{m = 0}^{4} \tilde{\omega}_{m} }, \;
	\text{with } \tilde{\omega}_{m} = \gamma _{m}(1 + \frac{\tau^{2}}{\beta_{m}^{2} + \varepsilon}), \;
	\text{for } m = 0, 1, 2, 3, 4,
	\label{eq:nonlinear-weights-2D}
\end{equation}
where $ \varepsilon = 10^{-8} $ is a small positive number to avoid division by zero.

Finally, a nonlinear HWENO reconstruction of the polynomial $ u_{i}(x, y) $ for $ u(x, y) $
within the target cell $ \tri_{i} $ is obtained by replacing some of the linear weights in
\cref{eq:rewrite-2D} by nonlinear weights defined in \cref{eq:nonlinear-weights-2D}
\begin{equation}
	u_{i,j}(x,y) = \omega_{0} (\frac{1}{\gamma _{0}}p_{0}(x,y) - \sum_{m = 1}^{4} \frac{\gamma_{m}}{\gamma _{0}}p_{m}(x,y))
	+ \sum_{m = 1}^{4} \omega_{m} p_{m}(x,y).
	\label{eq:hweno-reconstruction-2D}
\end{equation}
Those reconstructed point values at cell interfaces are evaluated by
\begin{equation}
	\left\{
	\begin{aligned}
		 & u_{i\pm\half,j}^{\mp,\mu} = u_{i,j}(x_{i\pm\half}, y_{j}^{\mu}),\quad u_{i,j\pm\half}^{\mu,\mp} = u_{i,j}(x_{i}^{\mu}, y_{j\pm\half}),                     \\
		 & v_{i\pm\half,j}^{\mp,\mu} = \pdv*{u_{i,j}}{x}(x_{i\pm\half}, y_{j}^{\mu}),\quad v_{i,j\pm\half}^{\mu,\mp} = \pdv*{u_{i,j}}{x}(x_{i}^{\mu}, y_{j\pm\half}), \\
		 & w_{i\pm\half,j}^{\mp,\mu} = \pdv*{u_{i,j}}{y}(x_{i\pm\half}, y_{j}^{\mu}),\quad w_{i,j\pm\half}^{\mu,\mp} = \pdv*{u_{i,j}}{y}(x_{i}^{\mu}, y_{j\pm\half}). \\
	\end{aligned}
	\right.
\end{equation}

Concurrently, the modified partial derivatives $ \hat{v}_{i} $ and $ \hat{w}_{i} $,
which will be only used in time discretization step,
are determined from $ u_{i}(x, y) $
by integrating its corresponding partial derivatives
over the target cell $ I_{i,j} $, we have
\begin{equation}
	\left\{
	\begin{aligned}
		 & \hat{v}_{i} = \frac{1}{\dy}
		\int_{I_{i,j}}^{} \pdv*{u_{i}(x,y)}{x}  \d{\Omega}, \quad \\
		 & \hat{w}_{i} = \frac{1}{\dx}
		\int_{I_{i,j}}^{} \pdv*{u_{i}(x,y)}{y}  \d{\Omega}. \quad \\
	\end{aligned}
	\right.
\end{equation}
\begin{step}
	Time discretization method for the semi-discrete scheme \cref{eq:2D-semi-discrete-scheme-uniform}.
\end{step}
As in the one-dimensional case, we employ the third-order strong stability preserving (SSP) Runge-Kutta
time discretization method to solve the semi-discrete scheme \cref{eq:2D-semi-discrete-scheme-uniform}.
Different from the one-dimensional case, the involved variables are $ \bar{u}_{i} $, $ \bar{v}_{i} $, and $ \bar{w}_{i} $.
And similarly, we modify the average partial derivatives $ \hat{v}_{i} $ and $ \hat{w}_{i} $
in each stage of the Runge-Kutta method.

\begin{equation}
	\left\{
	\begin{aligned}
		\vec{u}^{(1)}_{i,j}   & = \hat{\vec{u}}^{n}_{i,j} + \dt L_{i,j}(\vec{u}^{n}),                                                     \\
		\vec{u}^{(2)}_{i,j}   & = \frac{3}{4}\hat{\vec{u}}^{n}_{i,j} + \frac{1}{4}( \hat{\vec{u}}^{(1)}_{i,j} + \dt L_{i,j}(\vec{u}^{(1)}) ), \\
		\vec{u}^{n+1}_{i,j} & = \frac{1}{3}\hat{\vec{u}}^{n}_{i,j} + \frac{2}{3}( \hat{\vec{u}}^{(2)}_{i,j} + \dt L_{i,j}(\vec{u}^{(2)}) ), \\
	\end{aligned}
	\right.
\end{equation}
where $ \vec{u}^{\star}_{i,j} := [\bar{u}_{i,j}^{\star}, \bar{v}_{i,j}^{\star}, \bar{w}_{i,j}^{\star}] $,
$ \hat{\vec{u}}^{\star}_{i,j} := [\bar{u}_{i,j}^{\star}, \hat{v}_{i,j}^{\star}, \hat{w}_{i,j}^{\star}] $,
and $ L_{i,j} := [L_{i,j}^{u}, L_{i,j}^{v}, L_{i,j}^{w}] $.
These modified average partial derivatives are defined as
$ \hat{v}_{i,j}^{\star} = \frac{1}{\dx }\int_{I_{i,j}}^{} \pdv*{}{x} u(x,y,t^{\star}) \d{\Omega} $ and
$ \hat{w}_{i,j}^{\star} = \frac{1}{\dx }\int_{I_{i,j}}^{} \pdv*{}{y} u(x,y,t^{\star}) \d{\Omega} $,
where
$ \star $ denotes $ n, (1) $, and $ (2) $ corresponding to each time stage.

\begin{remark}
    For systems of conservation laws, such as the Euler equations of gas dynamics, the HWENO procedure of reconstruction and modification
	should be performed in the local characteristic directions to suppress spurious oscillations \cite{Shu1998EssentiallyNonoscillatory}, in both one and two dimensional cases.
\end{remark}

%% file: sections/sec-numerical-tests.tex
\section{Numerical Tests}
\label{sec:numerical-tests}
In this section, several benchmark examples are conducted to verify the high-order accuracy, efficiency, high resolution, and robustness of the constructed HWENO schemes.
Numerical results are compared with those of the fifth-order WENO schemes \cite{Zhu2017NewType, Zhu2017NewThird, Zhu2018NewFinite}
termed WENO, in both one- and two-dimensional cases.
The efficiency is compared with the first HWENO scheme of Qiu and Shu  \cite{Qiu2004HermiteWENOa}  termed HWENO-QS in the one-dimensional case.
For HWENO schemes, the linear weights for those linear polynomials are set as $ 1/400 $, in both one-dimensional and two-dimensional cases, and the quartic polynomial possesses the remaining weight. This choice follows the recommendation in \cite{Fan2025MomentBasedHermite}.
For fair comparison, the linear weights for WENO schemes are chosen as recommended in \cite{Zhu2017NewType, Zhu2018NewFinite},
meanwhile, the WENO schemes are implemented in the genuine two-dimensional reconstruction framework \cite{Zhu2018NewFinite}, instead of the dimension-by-dimension fashion \cite{Zhu2017NewType}.
The Courant-Friedrichs-Lewy (CFL) number is set as 0.6.

All numerical results were obtained by serial computations performed on a single physical core of an Intel Xeon Gold 6130 processor running at 2.10 GHz. The numerical code was written in Fortran 90 and compiled using the Intel Fortran Compiler (ifort), version 19.0.1.144,
with the \texttt{-O3} optimization option enabled.
\input{./sections/subsec-acc-test.tex}

\input{./sections/subsec-discontinuous-test.tex}

%% file: sections/subsec-acc-test.tex
\subsection{Accuracy tests}
In this subsection, we present several numerical examples with smooth solutions to verify the fifth-order accuracy of the constructed HWENO schemes.
	The numerical errors and computational costs of the HWENO schemes, the WENO schemes and the HWENO-QS scheme (in one-dimensional smooth case) are presented to reveal the better performance of the constructed HWENO schemes.

\begin{example}
	\label{ex:1d-burgers}
	We consider the one-dimensional nonlinear Burgers' equation with initial conditions:
	\begin{equation}
		\left\{
		\begin{aligned}
			 & u_{t} + (\frac{u^{2}}{2})_{x} = 0, \quad x \in [0,2], \\
			 & u(x,0) = 0.5 + \sin(\pi x), \quad                     \\
		\end{aligned}
		\right.
		\label{eq:1D-Burgers}
	\end{equation}
	where periodic boundary conditions are implemented.
	We simulate the problem up to the final time $ T = 0.5/\pi $, at which the solution remains smooth.
	We present the $ L^{1} $ and $ L^{\infty} $ errors and orders of accuracy of the constructed HWENO schemes and the WENO schemes.
	The results in \cref{tab:1d-burgers} indicate that all schemes achieve their designed fifth-order accuracies
	and the HWENO scheme produces less numerical errors with the same grid.
	The HWENO scheme is more efficient than the WENO scheme and the HWENO-QS scheme, as illustrated in \cref{fig:1d-burgers-efficiency}.
	\begin{table}[htbp]
\centering
\caption{\cref{ex:1d-burgers}. Burgers' equation, $L^1$, $L^\infty$ errors and orders, at $T=0.5/\pi$.}
\label{tab:1d-burgers}
\begin{tabular}{c c r r r r}
\toprule
Scheme & Mesh & $L^1$ error & Order & $L^\infty$ error & Order \\
\midrule
\multirow{6}{*}{HWENO} & $40$ & 2.66E-05 & -- & 1.39E-04 & -- \\
 & $80$ & 3.65E-07 & 6.19 & 4.32E-06 & 5.00 \\
 & $120$ & 4.64E-08 & 5.09 & 6.02E-07 & 4.86 \\
 & $160$ & 1.08E-08 & 5.06 & 1.44E-07 & 4.97 \\
 & $200$ & 3.57E-09 & 4.97 & 4.72E-08 & 5.01 \\
 & $240$ & 1.44E-09 & 5.00 & 1.90E-08 & 5.00 \\
\midrule
\multirow{6}{*}{WENO} & $40$ & 5.28E-05 & -- & 4.58E-04 & -- \\
 & $80$ & 1.89E-06 & 4.81 & 2.34E-05 & 4.29 \\
 & $120$ & 2.64E-07 & 4.85 & 3.41E-06 & 4.75 \\
 & $160$ & 6.35E-08 & 4.95 & 8.32E-07 & 4.90 \\
 & $200$ & 2.12E-08 & 4.92 & 2.75E-07 & 4.96 \\
 & $240$ & 8.55E-09 & 4.97 & 1.12E-07 & 4.93 \\
\midrule
\multirow{6}{*}{HWENO-QS} & $40$ & 3.79E-05 & -- & 3.48E-04 & -- \\
 & $80$ & 1.67E-06 & 4.50 & 1.57E-05 & 4.47 \\
 & $120$ & 2.60E-07 & 4.59 & 2.17E-06 & 4.87 \\
 & $160$ & 6.65E-08 & 4.74 & 6.51E-07 & 4.19 \\
 & $200$ & 2.24E-08 & 4.88 & 2.46E-07 & 4.37 \\
 & $240$ & 8.85E-09 & 5.08 & 1.08E-07 & 4.53 \\
\bottomrule
\end{tabular}
\end{table}
	\begin{figure}[htbp]
		\centering
		\includegraphics[width=0.495\textwidth]{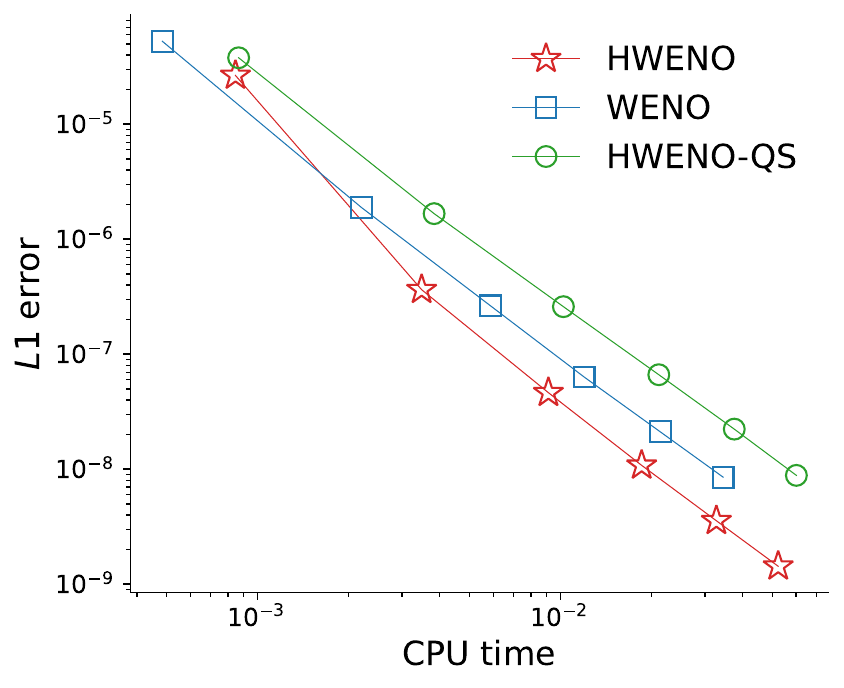}
		\includegraphics[width=0.495\textwidth]{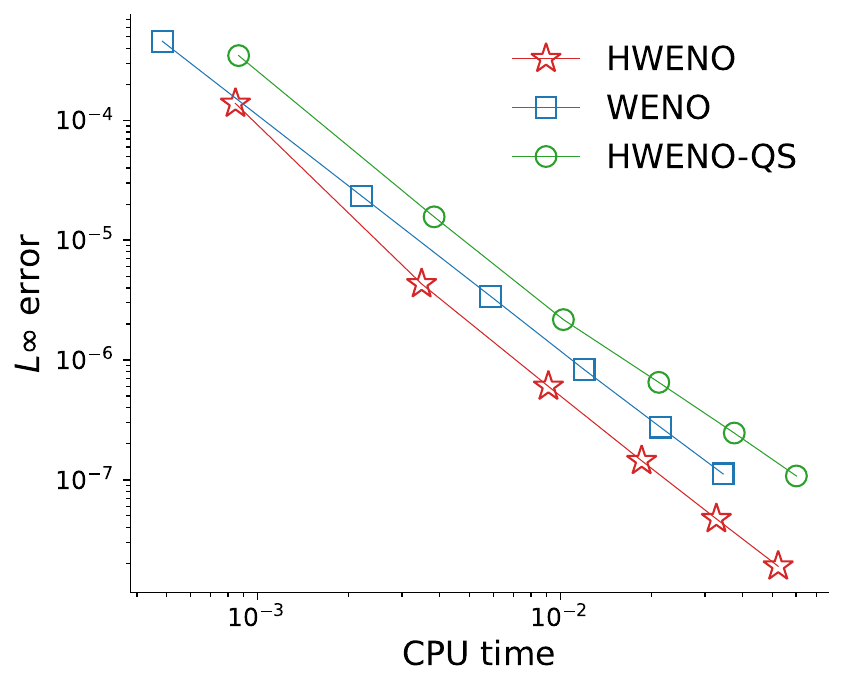}
		\caption{
		\cref{ex:1d-burgers}. Burgers' equation, comparison of $ L^{1} $ and $ L^{\infty} $ errors and CPU time.
	}
		\label{fig:1d-burgers-efficiency}
	\end{figure}
\end{example}
\begin{example}
	\label{ex:1d-euler}
	We solve the one-dimensional compressible Euler equations:
	\begin{equation}
		\pdv*{}{t}\begin{bmatrix}
			\rho \\ \rho \mu \\ E
		\end{bmatrix}
		+ \pdv*{}{x} \begin{bmatrix}
			\rho \mu \\ \rho \mu^{2} + p \\ \mu(E + p)
		\end{bmatrix}
		= 0
		\label{eq:1D-Euler}
	\end{equation}
	where $ \rho $, $ \mu $, $ E $, and $ p $ denote the density, velocity, total energy, and pressure, respectively.
	We solve the system in the computational domain $ [0, 2] $ with the initial conditions given by
	\begin{equation}
		(\rho, \mu, p, \gamma) = (1 + 0.2\sin(\pi x), 1, 1,1.4),
	\end{equation}
	and the periodic boundary conditions are applied.
	The exact solution at time $ t $ is $ (\rho, \mu, p) = (1+0.2\sin(\pi(x-t)), 1, 1) $.
	We simulate the system until time $ T = 2.0 $.
	The $ L^{1} $ and $ L^{\infty} $ errors and orders of accuracy of the proposed HWENO scheme
	and the classical WENO scheme are demonstrated in \cref{tab:1d-euler}.
	The results in \cref{tab:1d-euler} indicate that all schemes achieve their designed fifth-order accuracies
	and the HWENO scheme produces smaller numerical errors.
	The HWENO scheme is more efficient than the WENO scheme and the HWENO-QS scheme, as illustrated in \cref{fig:1d-euler-efficiency}.

	\begin{table}[htbp]
\centering
\caption{\cref{ex:1d-euler}. Euler equations, $L^1$, $L^\infty$ errors and orders, at $T=2.0$.}
\label{tab:1d-euler}
\begin{tabular}{c c r r r r}
\toprule
Scheme & Mesh & $L^1$ error & Order & $L^\infty$ error & Order \\
\midrule
\multirow{6}{*}{HWENO} & $40$ & 5.20E-07 & -- & 1.97E-06 & -- \\
 & $80$ & 1.62E-08 & 5.00 & 3.27E-08 & 5.91 \\
 & $120$ & 2.14E-09 & 5.00 & 3.74E-09 & 5.35 \\
 & $160$ & 5.07E-10 & 5.00 & 8.45E-10 & 5.17 \\
 & $200$ & 1.66E-10 & 5.00 & 2.73E-10 & 5.07 \\
 & $240$ & 6.65E-11 & 5.02 & 1.10E-10 & 4.98 \\
\midrule
\multirow{6}{*}{WENO} & $40$ & 3.09E-06 & -- & 5.25E-06 & -- \\
 & $80$ & 9.71E-08 & 4.99 & 1.55E-07 & 5.08 \\
 & $120$ & 1.28E-08 & 5.00 & 2.02E-08 & 5.02 \\
 & $160$ & 3.04E-09 & 5.00 & 4.79E-09 & 5.01 \\
 & $200$ & 9.95E-10 & 5.00 & 1.57E-09 & 5.00 \\
 & $240$ & 4.00E-10 & 5.00 & 6.32E-10 & 4.99 \\
\midrule
\multirow{6}{*}{HWENO-QS} & $40$ & 8.51E-06 & -- & 1.56E-05 & -- \\
 & $80$ & 2.61E-07 & 5.03 & 5.04E-07 & 4.95 \\
 & $120$ & 3.41E-08 & 5.02 & 6.46E-08 & 5.06 \\
 & $160$ & 8.02E-09 & 5.03 & 1.47E-08 & 5.14 \\
 & $200$ & 2.60E-09 & 5.04 & 4.58E-09 & 5.23 \\
 & $240$ & 1.03E-09 & 5.06 & 1.79E-09 & 5.14 \\
\bottomrule
\end{tabular}
\end{table}
	\begin{figure}[htbp]
		\centering
		\includegraphics[width=0.495\textwidth]{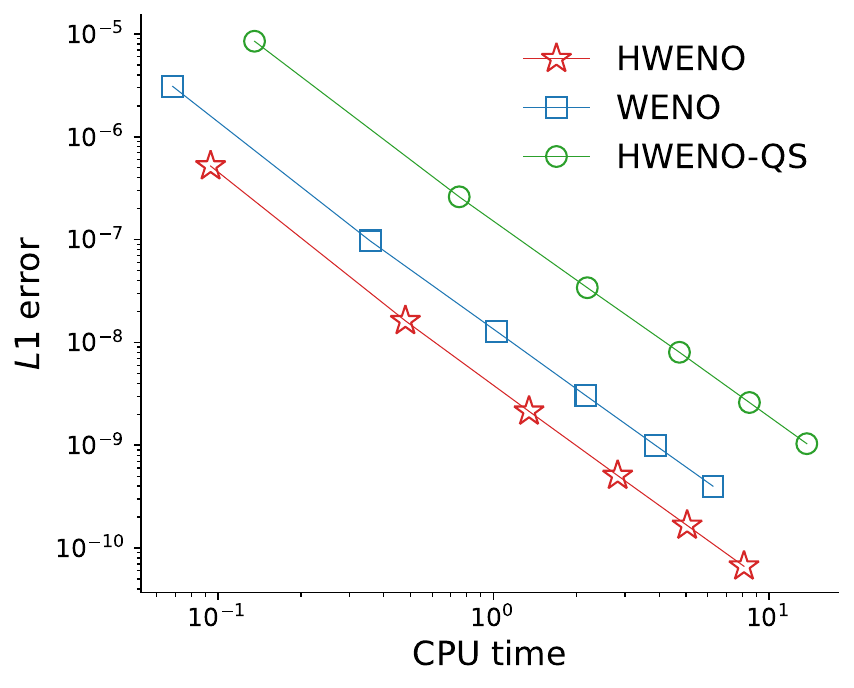}
		\includegraphics[width=0.495\textwidth]{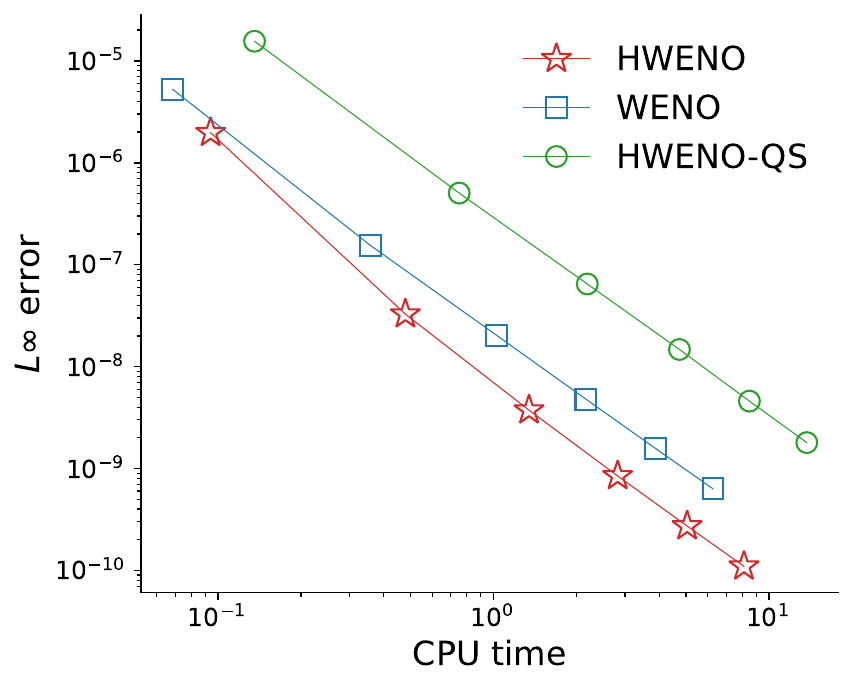}
		\caption{
		\cref{ex:1d-euler}. Euler equations, comparison of $ L^{1} $ and $ L^{\infty} $ errors and CPU time.
	}
		\label{fig:1d-euler-efficiency}
	\end{figure}
\end{example}
\begin{example}
	\label{ex:2d-burgers}
	We consider the two-dimensional nonlinear Burgers' equation with initial conditions:
	\begin{equation}
		\left\{
		\begin{aligned}
			 & u_{t} + (\frac{u^{2}}{2})_{x} + (\frac{u^{2}}{2})_{y} = 0, \quad (x, y) \in [0,4]\times [0,4], \\
			 & u(x,y,0) = 0.5 + \sin(0.5\pi(x + y)), \quad                                                    \\
		\end{aligned}
		\right.
		\label{eq:2D-Burgers}
	\end{equation}
	where periodic boundary conditions are implemented in each direction.
	We simulate this problem up to the final time $ T = 0.5/\pi $,
	at which the solution remains smooth.
	We present the $ L^{1} $ and $ L^{\infty} $ errors and orders of accuracy of the constructed HWENO schemes and the WENO schemes.
	The results in \cref{tab:2d-burgers} indicate that all schemes achieve their designed fifth-order accuracies
	and the HWENO scheme produces less numerical errors with the same grid.
	The HWENO scheme is more efficient than the WENO scheme, as illustrated in \cref{fig:2d-burgers-efficiency}.
	\begin{table}[htbp]
\centering
\caption{\cref{ex:2d-burgers}. Burgers' equation, $L^1$, $L^\infty$ errors and orders, at $T=0.5/\pi$.}
\label{tab:2d-burgers}
\begin{tabular}{c c r r r r}
\toprule
Scheme & Mesh & $L^1$ error & Order & $L^\infty$ error & Order \\
\midrule
\multirow{6}{*}{HWENO} & $40 \times 40$ & 3.87E-05 & -- & 2.51E-04 & -- \\
 & $80 \times 80$ & 1.13E-06 & 5.10 & 1.33E-05 & 4.24 \\
 & $120 \times 120$ & 1.52E-07 & 4.95 & 1.93E-06 & 4.77 \\
 & $160 \times 160$ & 3.71E-08 & 4.89 & 4.79E-07 & 4.84 \\
 & $200 \times 200$ & 1.23E-08 & 4.95 & 1.62E-07 & 4.87 \\
 & $240 \times 240$ & 4.98E-09 & 4.96 & 6.56E-08 & 4.94 \\
\midrule
\multirow{6}{*}{WENO} & $40 \times 40$ & 2.72E-04 & -- & 2.05E-03 & -- \\
 & $80 \times 80$ & 1.31E-05 & 4.38 & 1.47E-04 & 3.80 \\
 & $120 \times 120$ & 1.96E-06 & 4.68 & 2.38E-05 & 4.49 \\
 & $160 \times 160$ & 4.91E-07 & 4.81 & 6.21E-06 & 4.67 \\
 & $200 \times 200$ & 1.64E-07 & 4.91 & 2.12E-06 & 4.81 \\
 & $240 \times 240$ & 6.71E-08 & 4.91 & 8.69E-07 & 4.89 \\
\bottomrule
\end{tabular}
\end{table}
	\begin{figure}[htbp]
		\centering
		\includegraphics[width=0.495\textwidth]{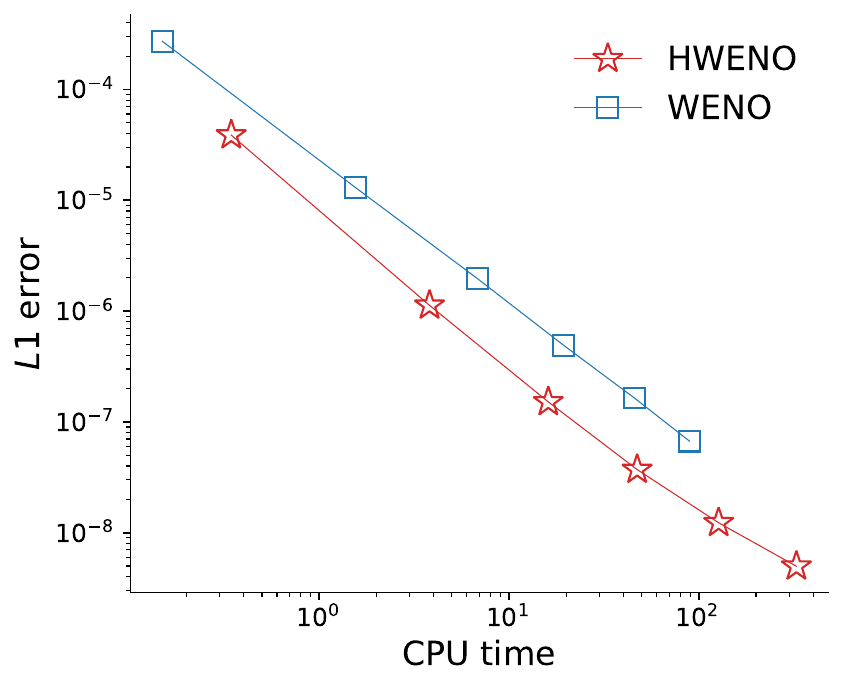}
		\includegraphics[width=0.495\textwidth]{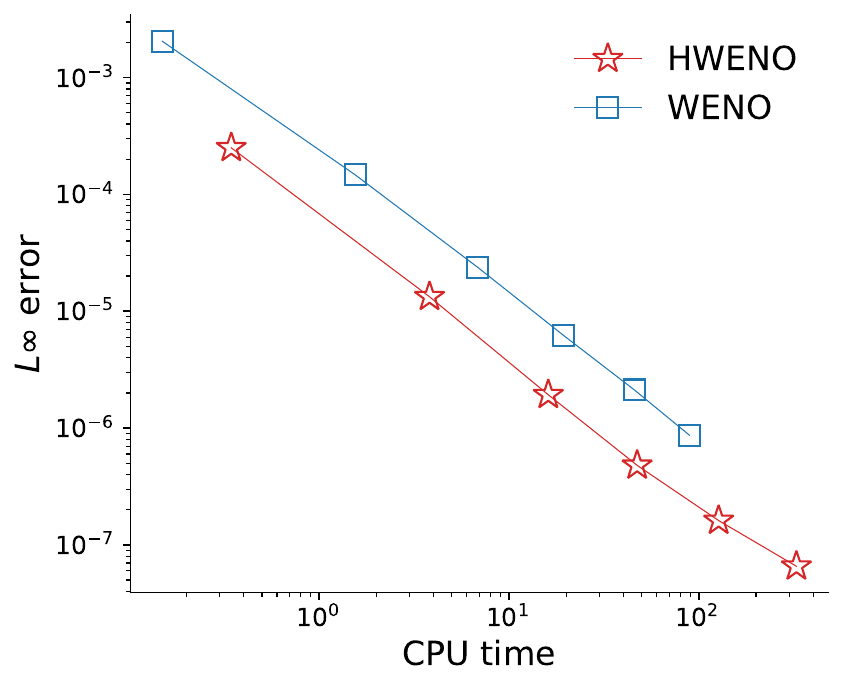}
		\caption{\cref{ex:2d-burgers}. Comparison of $ L^{1} $ and $ L^{\infty} $ errors and CPU time.}
		\label{fig:2d-burgers-efficiency}
	\end{figure}
\end{example}
\begin{example}
	\label{ex:2d-euler}
	We solve the two-dimensional compressible Euler equations:
	\begin{equation}
		\pdv*{}{t}\begin{bmatrix}
			\rho \\ \rho \mu \\ \rho \nu \\ E
		\end{bmatrix}
		+ \pdv*{}{x} \begin{bmatrix}
			\rho \mu \\ \rho \mu^{2} + p \\ \rho \mu \nu \\ \mu(E + p)
		\end{bmatrix}
		+ \pdv*{}{y} \begin{bmatrix}
			\rho \nu \\ \rho \mu \nu \\ \rho \nu^{2} + p \\ \nu(E + p)
		\end{bmatrix}
		= 0
		\label{eq:2D-Euler}
	\end{equation}
	where $ \rho $, $ \mu $, $ \nu $, $ E $, and $ p $ denote the density, the $x$-velocity, the $y$-velocity, the total energy, and the pressure, respectively.
	We solve the system in the computational domain $ [0, 2]\times[0, 2] $ with the initial conditions given by
	\begin{equation}
		(\rho, \mu, \nu, p, \gamma) = (1 + 0.2\sin(\pi(x + y)), 1, 1, 1, 1.4),
	\end{equation}
	and the periodic boundary conditions are applied in each direction.
	The exact solution at time $ t $ is
	$ (\rho, \mu, \nu, p) = (1 + 0.2\sin(\pi(x + y - 2t)), 1, 1, 1) $.
	We simulate the system until time $ T = 2.0 $.
	\begin{table}[htbp]
\centering
\caption{\cref{ex:2d-euler}. Euler equations, $L^1$, $L^\infty$ errors and orders, at $T=2.0$.}
\label{tab:2d-euler}
\begin{tabular}{c c r r r r}
\toprule
Scheme & Mesh & $L^1$ error & Order & $L^\infty$ error & Order \\
\midrule
\multirow{6}{*}{HWENO} & $40 \times 40$ & 3.61E-06 & -- & 7.12E-06 & -- \\
 & $80 \times 80$ & 1.13E-07 & 4.99 & 1.86E-07 & 5.26 \\
 & $120 \times 120$ & 1.49E-08 & 5.00 & 2.38E-08 & 5.07 \\
 & $160 \times 160$ & 3.54E-09 & 5.00 & 5.61E-09 & 5.03 \\
 & $200 \times 200$ & 1.16E-09 & 5.00 & 1.84E-09 & 5.00 \\
 & $240 \times 240$ & 4.67E-10 & 4.99 & 7.44E-10 & 4.96 \\
\midrule
\multirow{6}{*}{WENO} & $40 \times 40$ & 4.89E-05 & -- & 7.91E-05 & -- \\
 & $80 \times 80$ & 1.55E-06 & 4.98 & 2.44E-06 & 5.02 \\
 & $120 \times 120$ & 2.04E-07 & 4.99 & 3.22E-07 & 5.00 \\
 & $160 \times 160$ & 4.85E-08 & 5.00 & 7.63E-08 & 5.00 \\
 & $200 \times 200$ & 1.59E-08 & 5.00 & 2.50E-08 & 5.00 \\
 & $240 \times 240$ & 6.40E-09 & 5.00 & 1.01E-08 & 5.00 \\
\bottomrule
\end{tabular}
\end{table}
	\begin{figure}[htbp]
		\centering
		\includegraphics[width=0.495\textwidth]{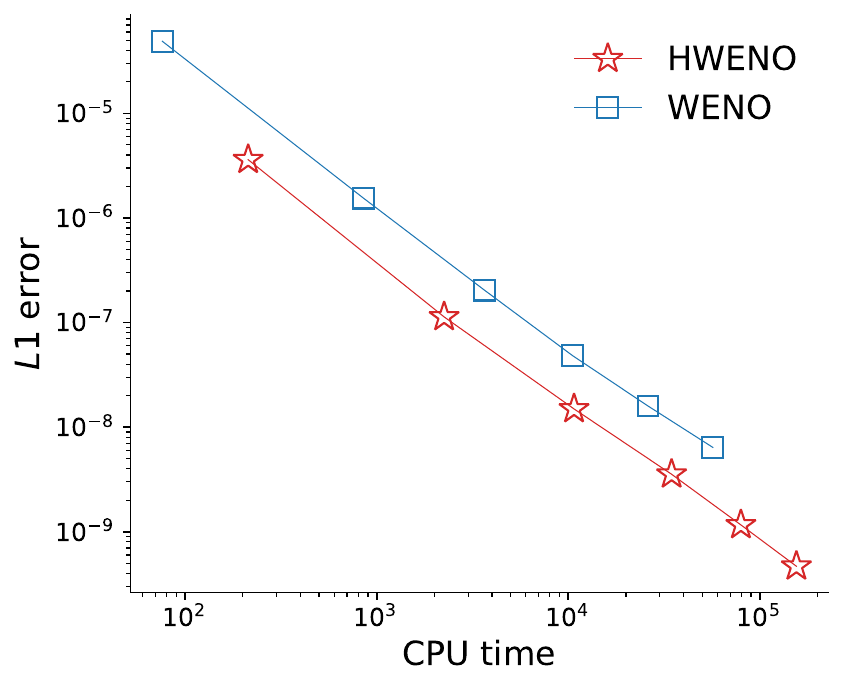}
		\includegraphics[width=0.495\textwidth]{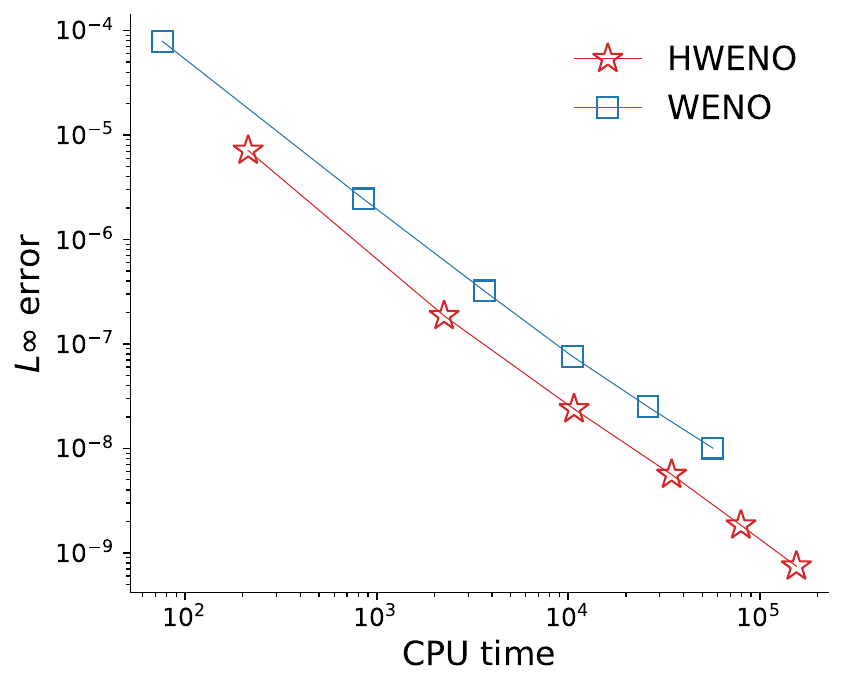}
		\caption{\cref{ex:2d-euler}. Euler equations, comparison of $ L^{1} $ and $ L^{\infty} $ errors and CPU time.}
		\label{fig:triangular-euler-efficiency}
	\end{figure}
\end{example}

%% file: sections/subsec-discontinuous-test.tex
\subsection{Discontinuous tests}
In this subsection, several benchmark examples involving discontinuities are conducted to verify the high resolution and the ability to capture shocks of the constructed HWENO schemes. It is also illustrated that the HWENO schemes behave better than the WENO schemes when simulations contain shocks and discontinuities.
\begin{example}
	\label{ex:1D-lax}
	We solve the one-dimensional Lax problem modeled by the Euler equations \cref{eq:1D-Euler}
	with initial conditions:
	\begin{equation}
		(\rho, \mu, p,\gamma) = \left\{
		\begin{aligned}
			 & (0.445,0.698,3.528,1.4), \quad x\in [-0.5,0], \\
			 & (0.5,0,0.571,1.4), \quad       x\in [0,0.5].  \\
		\end{aligned}
		\right.
	\end{equation}
	The computational domain is $ [-0.5, 0.5] $ and the final time is $ T = 0.16 $.
	The inflow and outflow boundary conditions are imposed on the left and right boundaries, respectively.
	The HWENO scheme achieves higher resolution than the WENO scheme near discontinuities, as depicted in \cref{fig:lax}.
\begin{figure}[htbp]
	\centering
	\begin{subfigure}{.495\textwidth}
		\centering
		\includegraphics[width=\linewidth]{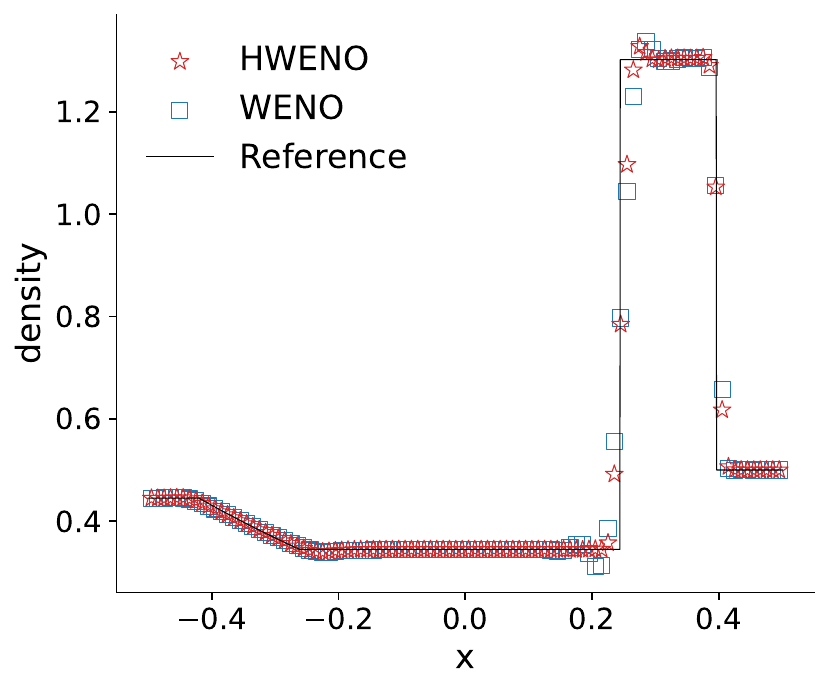}
	\end{subfigure}
	\begin{subfigure}{.495\textwidth}
		\centering
		\includegraphics[width=\linewidth]{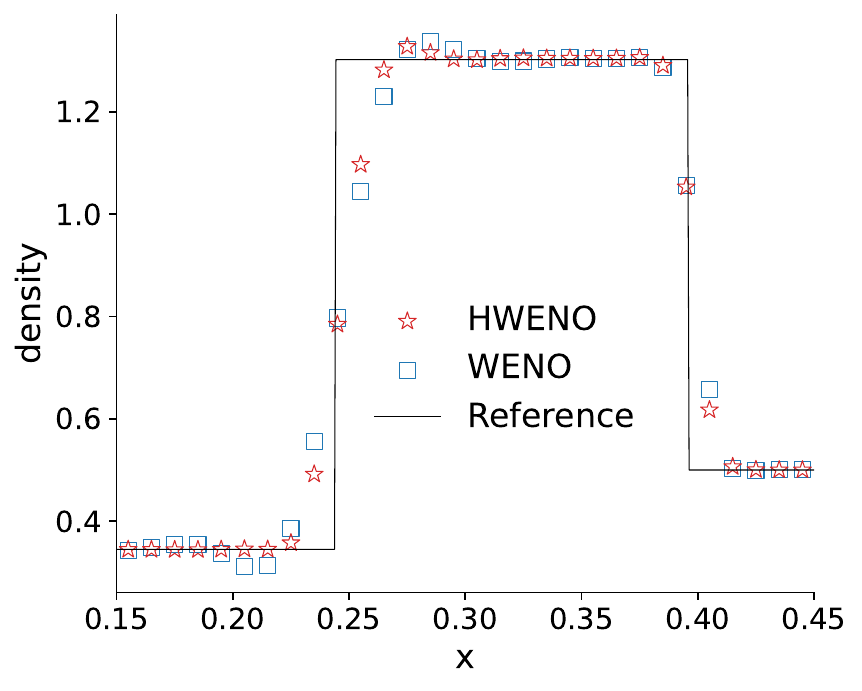}
	\end{subfigure}
	\caption{
	\cref{ex:1D-lax}. Lax problem, computed density of HWENO and WENO schemes with $ N = 100 $ cells at $ T = 0.16 $.
}
	\label{fig:lax}
\end{figure}
\end{example}
\begin{example}
	\label{ex:1D-blast}
	We solve the one-dimensional blast wave problem modeled by the Euler equations \cref{eq:1D-Euler}
	with initial conditions:
	\begin{equation}
		(\rho, \mu, p,\gamma) = \left\{
		\begin{aligned}
			 & (1,0,1000,1.4), \quad x\in [0,0.1],   \\
			 & (1,0,0.01,1.4), \quad x\in [0.1,0.9], \\
			 & (1,0,100,1.4), \quad x\in [0.9,1].    \\
		\end{aligned}
		\right.
	\end{equation}
	The computational domain is $ [0,1] $ and the final time is $ T = 0.038 $.
	The reflective boundary condition is imposed on all boundaries.
	The HWENO scheme achieves higher resolution than the WENO scheme near discontinuities, as depicted in \cref{fig:blast}.
\begin{figure}[htbp]
	\centering
	\begin{subfigure}{.495\textwidth}
		\centering
		\includegraphics[width=\linewidth]{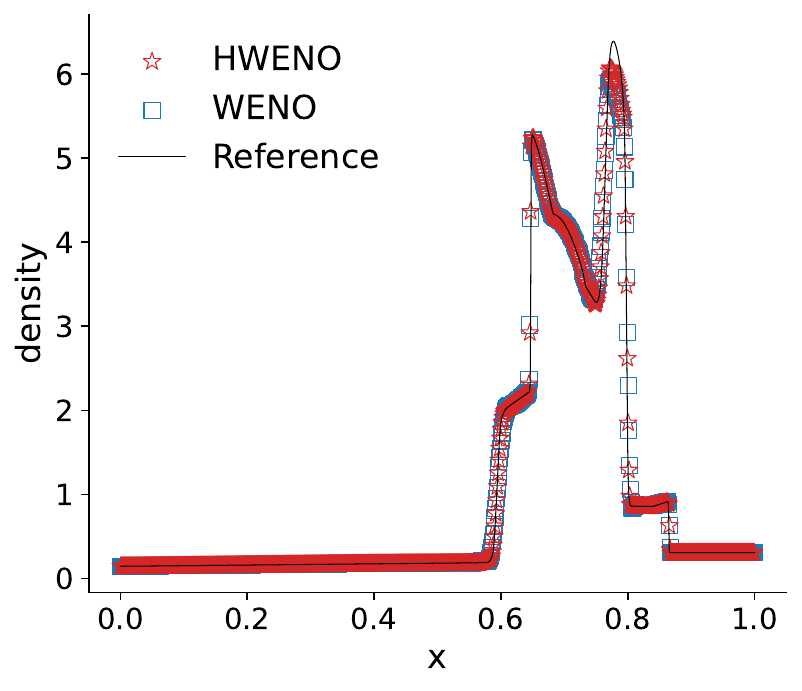}
	\end{subfigure}
	\begin{subfigure}{.495\textwidth}
		\centering
		\includegraphics[width=\linewidth]{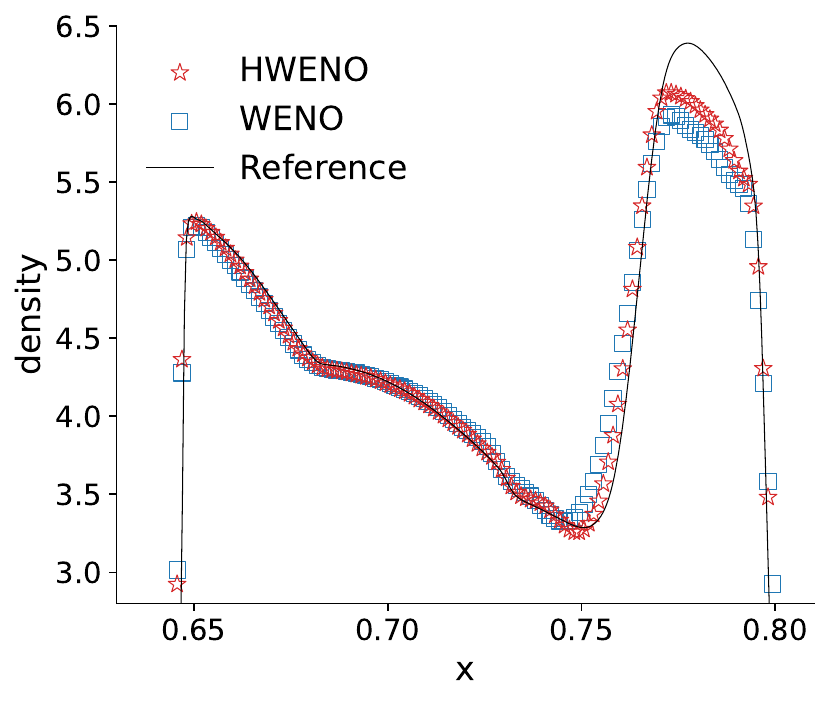}
	\end{subfigure}
	\caption{
		\cref{ex:1D-blast}. Blast problem, computed density of HWENO and WENO schemes with $ N = 800 $ cells at $ T = 0.038 $.
	}
	\label{fig:blast}
\end{figure}
\end{example}
\begin{example}
	\label{ex:1D-shu-osher}
	We solve the one-dimensional Shu-Osher problem modeled by the Euler equations \cref{eq:1D-Euler}
	with initial conditions:
	\begin{equation}
		(\rho, \mu, p,\gamma) = \left\{
		\begin{aligned}
			 & (3.857143,2.629369,10.33333,1.4), \quad x\in [-5,-4], \\
			 & (1 + 0.2\sin(5\pi x),0,1,1.4), \quad x\in [-4,5].   \\
		\end{aligned}
		\right.
	\end{equation}
	The computational domain is set as $ [-5,5] $ and the final time is $ T = 1.8 $.
	The outflow boundary condition is imposed on all boundaries.
	The HWENO scheme is comparable to the WENO scheme, as depicted in \cref{fig:shu-osher}.
\begin{figure}[htbp]
	\centering
	\begin{subfigure}{.495\textwidth}
		\centering
		\includegraphics[width=\linewidth]{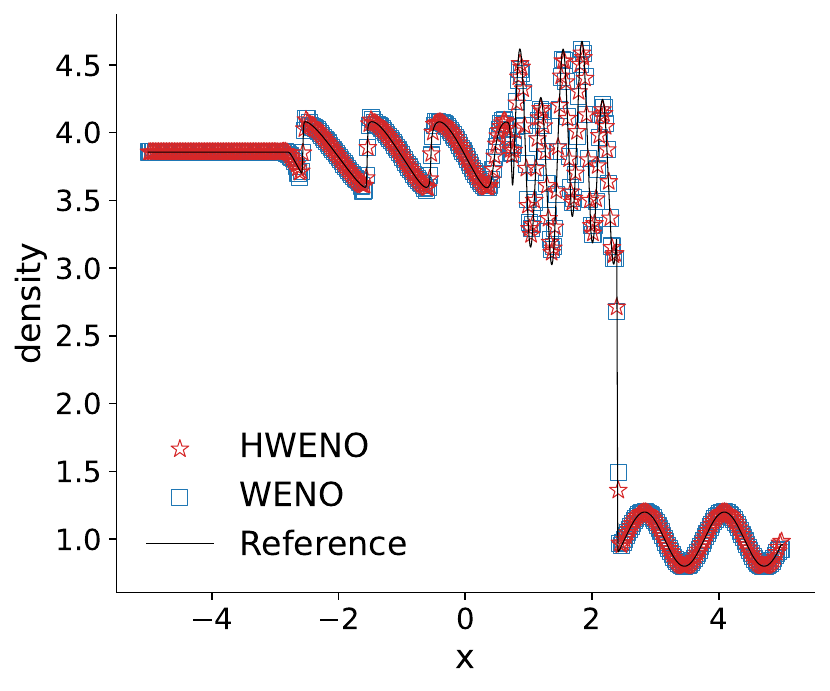}
	\end{subfigure}
	\begin{subfigure}{.495\textwidth}
		\centering
		\includegraphics[width=\linewidth]{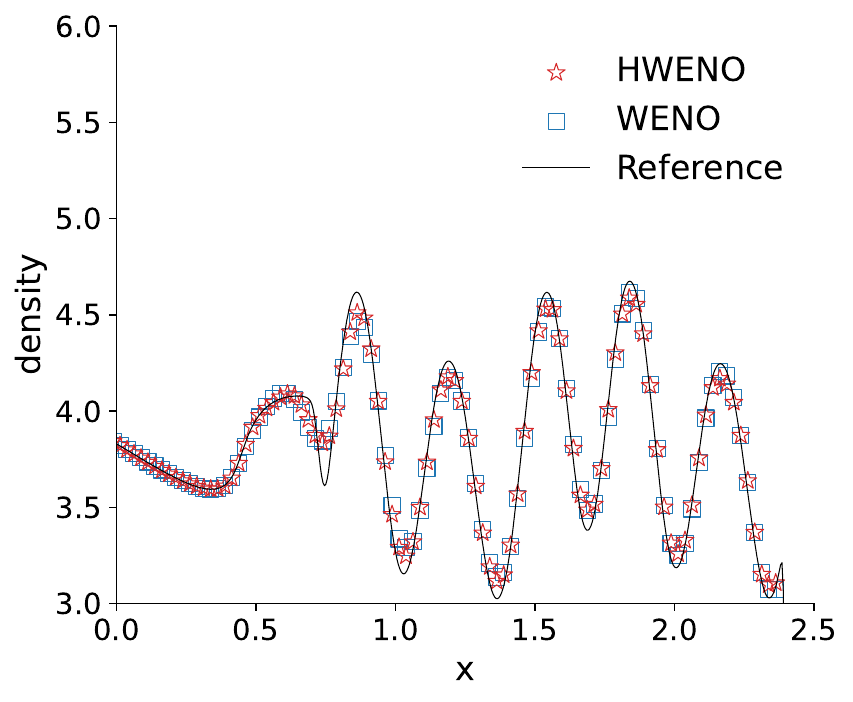}
	\end{subfigure}
	\caption{
	\cref{ex:1D-shu-osher}. Shu-Osher problem, computed density of HWENO and WENO schemes with $ N = 400 $ cells at $ T = 1.8 $.
}
	\label{fig:shu-osher}
\end{figure}
\end{example}
\begin{example}
	\label{ex:double-mach}
	We solve the double Mach reflection problem \cite{Woodward1984NumericalSimulation}, modeled by the Euler equations \cref{eq:2D-Euler}.
	The computational domain is set as $ [0,4]\times[0,1] $
	and the initial condition is
	\begin{equation}
		(\rho, \mu, \nu, p, \gamma) = \left\{
		\begin{aligned}
			 & (8,\frac{33}{4}\sin\frac{\pi}{3},-\frac{33}{4}\cos\frac{\pi}{3},116.5,1.4), \quad
			 x < \frac{1}{6} + \frac{y}{\sqrt{3}}, \\
			 & (1.4,0,0,1,1.4), \quad          \text{otherwise.}         \\
		\end{aligned}
		\right.
	\end{equation}
	Inflow and outflow boundary conditions are applied to the left and right boundaries, respectively.
	An exact post-shock boundary condition is imposed on the bottom from $ x = 0 $ to $ x = 1/6 $,
	and the remaining bottom portion is subjected to a reflective boundary condition.
	Additionally, the motion of a Mach 10 shock is applied to the top boundary.
	The simulation stops at a final time of $ T = 0.2 $.
	As depicted in \cref{fig:double-mach}, the HWENO scheme captures more complicated shock structures than the WENO scheme.
\begin{figure}[H]
	\centering
	\begin{subfigure}{\textwidth}
		\centering
	\includegraphics[width=\textwidth]{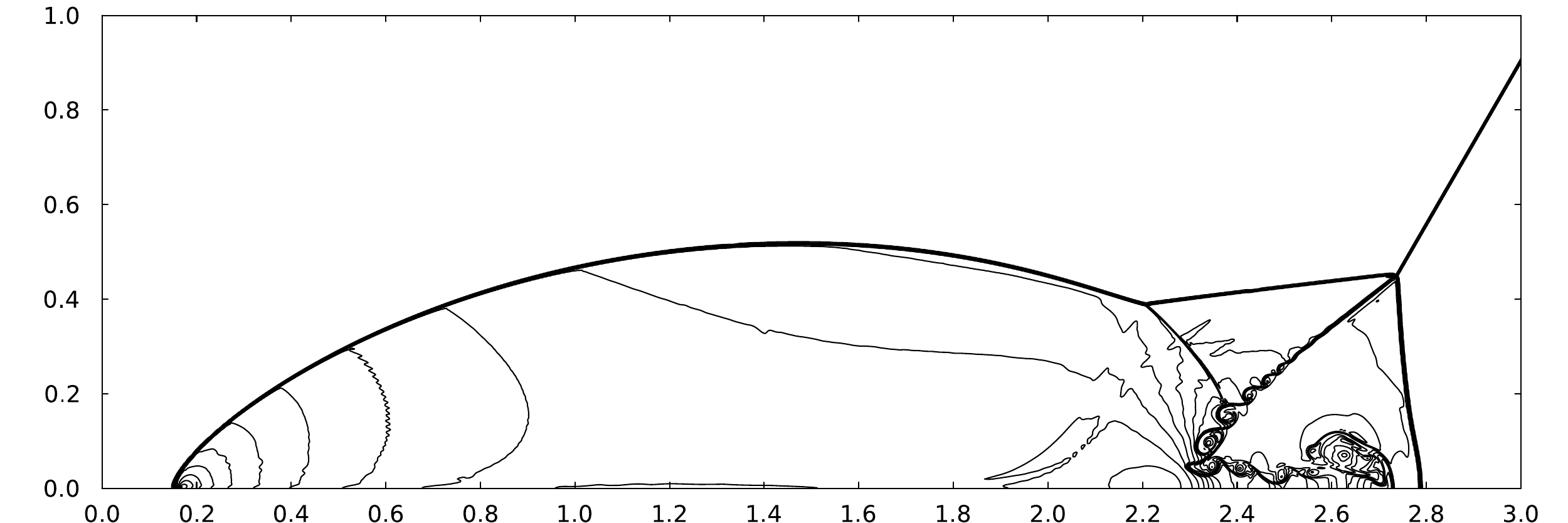}
		\caption{HWENO}
	\end{subfigure}

	\begin{subfigure}{\textwidth}
		\centering
	\includegraphics[width=\textwidth]{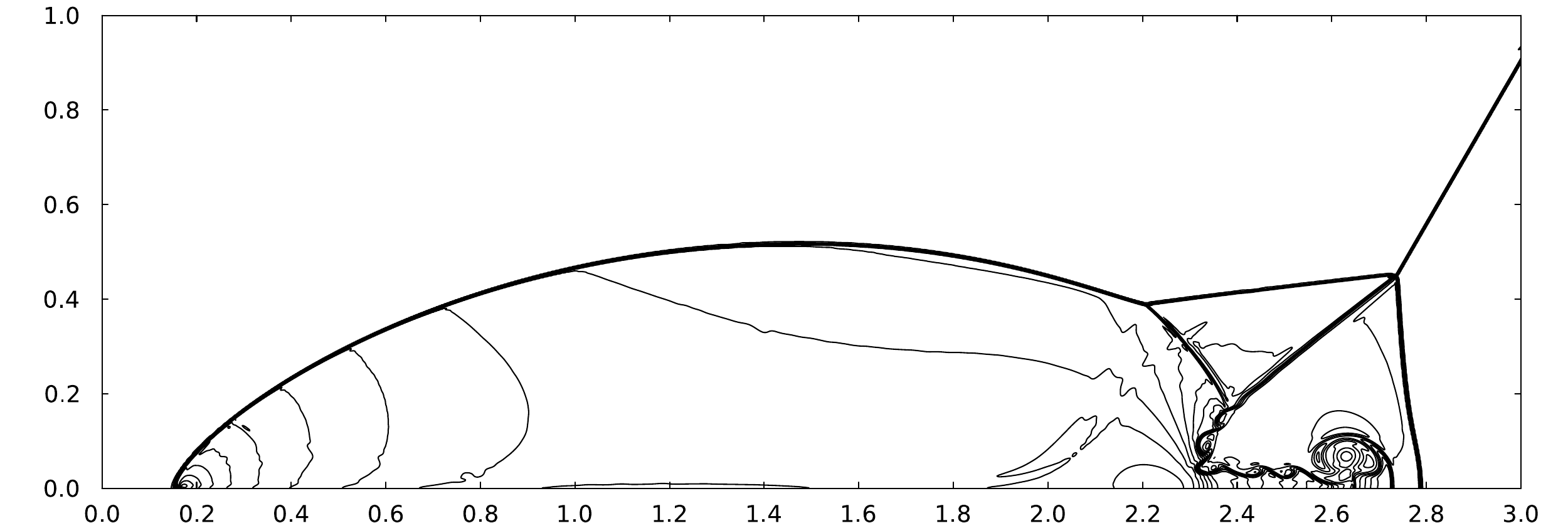}
		\caption{WENO}
	\end{subfigure}
	\caption{\cref{ex:double-mach}. Double Mach reflection problem, density contours with 30 equally spaced levels from 1.5 to 22.7. Mesh size: $ 1920\times 480 $.}
	\label{fig:double-mach}
\end{figure}
\end{example}
\begin{example}
	\label{ex:forward-step}
	We finally address the forward step problem \cite{Woodward1984NumericalSimulation}, modeled by the Euler equations \cref{eq:2D-Euler}
	on the computational domain of $ [0,0.6]\times[0,1] \cup [0.6, 1]\times[0.2,1] $,
	featuring a Mach 3 wind tunnel incorporating a step.
	The initial condition consists of a right-propagating Mach 3 flow,
	and the final time is set as $ T = 4 $.
	Reflective boundary conditions are imposed on the tunnel walls and the step,
	while inflow and outflow boundary conditions are specified on the left and right boundaries, respectively. For singularity at the corner of the step, we follow the treatment in \cite{Woodward1984NumericalSimulation}.
	As shown in \cref{fig:forward-step-coarse}, the HWENO scheme achieves better resolution than the WENO scheme.
	In the multi-resolution HWENO scheme proposed by \textcite{Li2021MultiresolutionHWENO},
	the unified set of stencils is also employed to construct
	a derivative-based sixth-order multi-resolution HWENO scheme in the finite volume framework and a fourth-order one in the finite difference framework.
	And they mentioned that their sixth-order finite volume multi-resolution HWENO scheme does not behave better than
	the fourth-order finite difference one in the simulation of this example.
	Our fifth-order finite volume HWENO scheme does not have this issue as shown in \cref{fig:forward-step-coarse}.

\begin{figure}[htbp]
	\centering
	\begin{subfigure}{\textwidth}
		\centering
	\includegraphics[width=\textwidth]{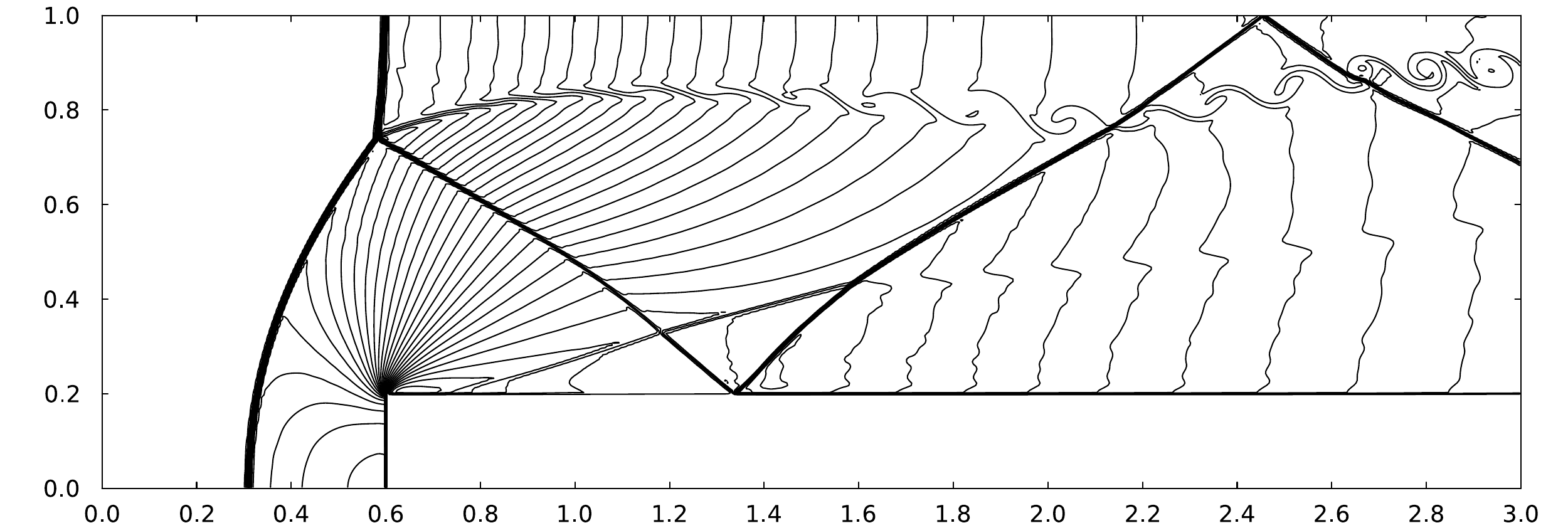}
		\caption{HWENO}
	\end{subfigure}

	\begin{subfigure}{\textwidth}
		\centering
	\includegraphics[width=\textwidth]{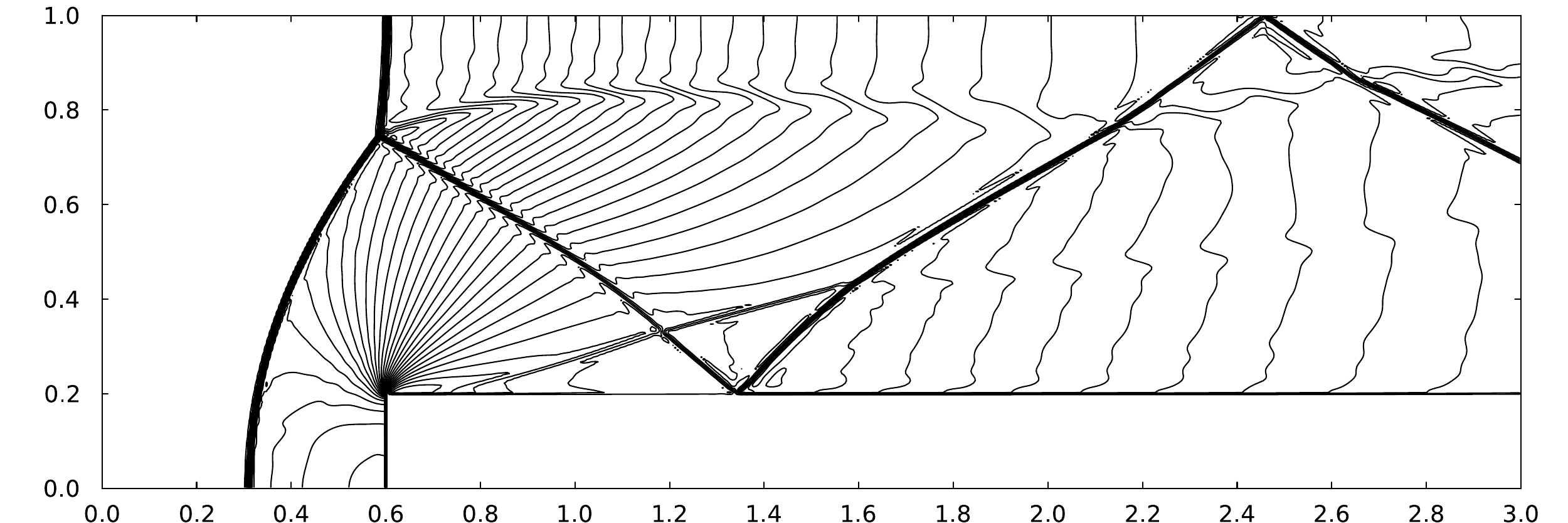}
		\caption{WENO}
	\end{subfigure}

	\begin{subfigure}{\textwidth}
		\centering
	\includegraphics[width=\textwidth]{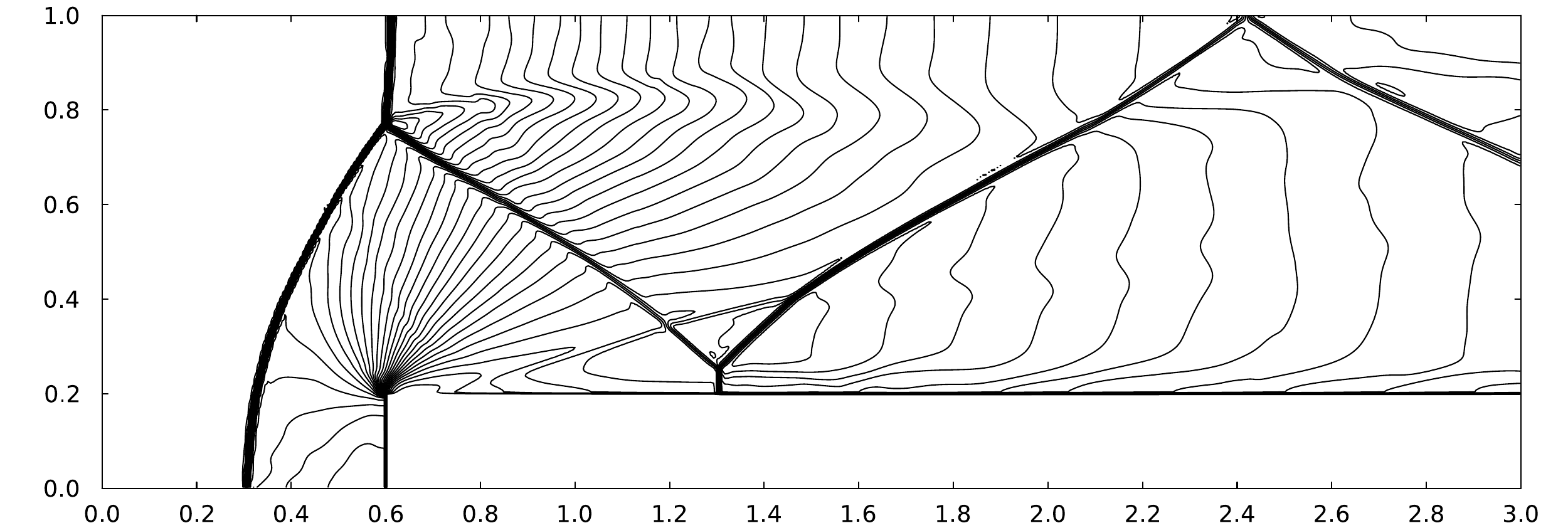}
	\caption{sixth-order multi-resolution HWENO \cite{Li2021MultiresolutionHWENO} }
	\end{subfigure}
	\caption{
	\cref{ex:forward-step}. Forward step problem, density contours with 30 equally spaced levels from 0.32 to 6.15. Mesh size: $ 600\times 200$.
}
	\label{fig:forward-step-coarse}
\end{figure}
\end{example}

%% file: sections/sec-conclusion.tex
\section{Concluding Remarks}
\label{sec:conclusions}
In this paper, a fifth-order derivative-based HWENO scheme is constructed for solving hyperbolic conservation laws.
The constructed scheme retains the merits of the previous HWENO schemes, such as the compact and unified stencils and artificial linear weights, and meanwhile it is equipped with an efficient and effective smoothness indicator.
The essential ingredient of the proposed scheme is that the information of the target cell's derivatives is excluded from the spatial reconstruction, while the same reconstructed polynomial is employed both in the spatial discretization and in modifying the derivatives during the temporal evolution.
Consequently, a single unified set of reconstruction stencils is used throughout the whole algorithm, which simplifies the implementation, reduces the computational cost, and improves the robustness near discontinuities.
Moreover, the scheme allows arbitrary positive linear weights that sum to one, while retaining a more compact stencil than the WENO scheme of the same order.

Different from the first derivative-based finite volume HWENO scheme \cite{Qiu2004HermiteWENOa}, which employs two different sets of reconstruction stencils for the governing equation and the derivative equation, the proposed scheme uses a unified stencil.
Consequently, it is not only easier to implement and more efficient, but also more robust near discontinuities.
Compared with the multi-resolution HWENO scheme \cite{Li2021MultiresolutionHWENO}, which also adopts the unified stencil but suffers from insufficient control of derivative values near strong shock waves, the proposed scheme achieves better resolution and robustness in such problems with strong shocks.
Compared with the WENO scheme \cite{Zhu2018NewFinite}, the proposed HWENO scheme achieves better efficiency and resolution.
Extensive numerical experiments confirm that the scheme achieves the designed fifth-order accuracy, produces smaller errors with less CPU time than the compared WENO scheme \cite{Zhu2018NewFinite} and
the first HWENO scheme \cite{Qiu2004HermiteWENOa},
and captures strong discontinuities with higher resolution and better robustness than the multi-resolution HWENO scheme \cite{Li2021MultiresolutionHWENO}.

%% file: sections/sec-acknowledgements.tex
\section*{CRediT authorship contribution statement}

\textbf{Peiwen Chen}: Writing -- review \& editing, Writing -- original draft, Software, Methodology, Investigation, Conceptualization.
\textbf{Zhuang Zhao}: Writing -- review \& editing, Methodology, Investigation, Conceptualization, Supervision, Funding acquisition.

\section*{Declaration of competing interest}

The authors declare that they have no known competing financial interests or personal relationships that could have appeared to influence the work reported in this paper.

\section*{Acknowledgements}

The research was partially supported by National Key R\&D Program of China under Grant Number 2024YFA1012500, National Natural Science Foundation of China
under Grant Number 12401541, and Fujian Provincial Natural Science Foundation of China under Grant Number 2026J009009.

\section*{Data availability}

No data were used for the research described in the article.

%% file: sections/sec-coefficients.tex
\section{The Coefficients of the Reconstructed Polynomials}
\label{sec:coefficients}
In the one-dimensional case, the coefficients of the reconstructed polynomials $ \{p_{m}(x)\}_{m=0}^{2} $ in \cref{eq:hweno-reconstruction-1D} are given as follows:
\begin{equation}
	\left\{
	\begin{aligned}
		c_{0,0} & = -\frac{47 }{480}\bar{u}_{i-1}+\frac{287 }{240}\bar{u}_i-\frac{47 }{480}\bar{u}_{i+1} -\frac{9\dx}{320}  \bar{v}_{i-1}+\frac{9\dx}{320}  \bar{v}_{i+1}, \\
		c_{0,1} & = -\frac{13}{16}  \bar{u}_{i-1}+\frac{13 }{16}\bar{u}_{i+1}-\frac{5\dx}{16}  \bar{v}_{i-1}-\frac{5\dx}{16}  \bar{v}_{i+1},                               \\
		c_{0,2} & =\frac{5 }{4}\bar{u}_{i-1}-\frac{5 }{2}\bar{u}_i+\frac{5 }{4}\bar{u}_{i+1}+ \frac{3\dx}{8}  \bar{v}_{i-1}-\frac{3\dx}{8}  \bar{v}_{i+1},                 \\
		c_{0,3} & =\frac{1}{4}\bar{u}_{i-1}-\frac{1}{4}\bar{u}_{i+1}+\frac{\dx}{4}  \bar{v}_{i-1}+\frac{\dx}{4}  \bar{v}_{i+1},                                            \\
		c_{0,4} & =-\frac{1}{2}\bar{u}_{i-1}+\bar{u}_i-\frac{1}{2}\bar{u}_{i+1}-\frac{\dx}{4}  \bar{v}_{i-1}+\frac{\dx}{4}  \bar{v}_{i+1},                                 \\
		c_{1,0} & = \bar{u}_{i},\quad c_{1,1} = \bar{u}_{i}-\bar{u}_{i-1},                                                                                                 \\
		c_{2,0} & = \bar{u}_{i},\quad c_{2,1} = \bar{u}_{i+1}-\bar{u}_{i},
	\end{aligned}
	\right.
\end{equation}

In the 2-dimensional case, the coefficients of the reconstructed polynomials $ \{p_{m}(x, y)\}_{m=0}^{4} $ in \cref{eq:hweno-reconstruction-2D} are given as follows:

\begin{equation}
	\left\{
	\begin{aligned}
		c_{0,1}  & =   (5\bar{q}_1-292\bar{q}_2+5\bar{q}_3-292\bar{q}_4+4028\bar{q}_5-292\bar{q}_6                                                                      \\
		         & +5\bar{q}_7-292\bar{q}_8+5\bar{q}_9-81\bar{q}_{10}+81\bar{q}_{11}-81\bar{q}_{16}+81\bar{q}_{17})\,\big/\,2880,                                       \\
		c_{0,2}  & =   \frac{1}{24}(2\bar{q}_1-2\bar{q}_3-19\bar{q}_4+19\bar{q}_6+2\bar{q}_7-2\bar{q}_9-6\bar{q}_{10}-6\bar{q}_{11}+3\bar{q}_{12}+3\bar{q}_{13}),       \\
		c_{0,3}  & =   \frac{1}{24}(2\bar{q}_1-19\bar{q}_2+2\bar{q}_3-2\bar{q}_7+19\bar{q}_8-2\bar{q}_9+3\bar{q}_{14}+3\bar{q}_{15}-6\bar{q}_{16}-6\bar{q}_{17}),       \\
		c_{0,4}  & =   \frac{1}{48}(-\bar{q}_1+2\bar{q}_2-\bar{q}_3+62\bar{q}_4-124\bar{q}_5+62\bar{q}_6-\bar{q}_7+2\bar{q}_8-\bar{q}_9+18\bar{q}_{10}-18\bar{q}_{11}), \\
		c_{0,5}  & =   \frac{1}{8}(-3\bar{q}_1+3\bar{q}_3+3\bar{q}_7-3\bar{q}_9-5\bar{q}_{12}+5\bar{q}_{13}-5\bar{q}_{14}+5\bar{q}_{15}),                               \\
		c_{0,6}  & =   \frac{1}{48}(-\bar{q}_1+62\bar{q}_2-\bar{q}_3+2\bar{q}_4-124\bar{q}_5+2\bar{q}_6-\bar{q}_7+62\bar{q}_8-\bar{q}_9+18\bar{q}_{16}-18\bar{q}_{17}), \\
		c_{0,7}  & =   \frac{1}{20}(-\bar{q}_1+\bar{q}_3+4\bar{q}_4-4\bar{q}_6-\bar{q}_7+\bar{q}_9+4\bar{q}_{10}+4\bar{q}_{11}-2\bar{q}_{12}-2\bar{q}_{13}),            \\
		c_{0,8}  & =   \frac{1}{4}(-\bar{q}_1+2\bar{q}_2-\bar{q}_3+\bar{q}_7-2\bar{q}_8+\bar{q}_9),                                                                     \\
		c_{0,9}  & =   \frac{1}{4}(-\bar{q}_1+\bar{q}_3+2\bar{q}_4-2\bar{q}_6-\bar{q}_7+\bar{q}_9),                                                                     \\
		c_{0,10} & =  \frac{1}{20}(-\bar{q}_1+4\bar{q}_2-\bar{q}_3+\bar{q}_7-4\bar{q}_8+\bar{q}_9-2\bar{q}_{14}-2\bar{q}_{15}+4\bar{q}_{16}+4\bar{q}_{17}),             \\
		c_{0,11} & =  \frac{1}{4}(-2\bar{q}_4+4\bar{q}_5-2\bar{q}_6-\bar{q}_{10}+\bar{q}_{11}),                                                                         \\
		c_{0,12} & =  \frac{1}{4}(\bar{q}_1-\bar{q}_3-\bar{q}_7+\bar{q}_9+2\bar{q}_{12}-2\bar{q}_{13}),                                                                 \\
		c_{0,13} & =  \frac{1}{4}(\bar{q}_1-2\bar{q}_2+\bar{q}_3-2\bar{q}_4+4\bar{q}_5-2\bar{q}_6+\bar{q}_7-2\bar{q}_8+\bar{q}_9),                                      \\
		c_{0,14} & =  \frac{1}{4}(\bar{q}_1-\bar{q}_3-\bar{q}_7+\bar{q}_9+2\bar{q}_{14}-2\bar{q}_{15}),                                                                 \\
		c_{0,15} & =  \frac{1}{4}(-2\bar{q}_2+4\bar{q}_5-2\bar{q}_8-\bar{q}_{16}+\bar{q}_{17}),                                                                         \\
		c_{1,1}  & =  \bar{u}_{i,j},\quad c_{1,2} = \bar{u}_{i,j}-\bar{u}_{i-1,j},\quad c_{1,3} = \bar{u}_{i,j} - \bar{u}_{i,j-1}                                       \\
		c_{2,1}  & =  \bar{u}_{i,j},\quad c_{2,2} = \bar{u}_{i+1,j}-\bar{u}_{i,j},\quad c_{2,3} =  \bar{u}_{i,j} - \bar{u}_{i,j-1}                                      \\
		c_{3,1}  & =  \bar{u}_{i,j},\quad c_{3,2} = \bar{u}_{i+1,j}-\bar{u}_{i,j},\quad c_{3,3} = \bar{u}_{i,j+1}-\bar{u}_{i,j}                                                       \\
		c_{4,1}  & =  \bar{u}_{i,j},\quad c_{4,2} = \bar{u}_{i,j}-\bar{u}_{i-1,j},\quad c_{4,3} = \bar{u}_{i,j+1}-\bar{u}_{i,j}
	\end{aligned}
	\right.
\end{equation}
where $ \bar{q}_{1}, \bar{q}_{2},\ldots, \bar{q}_{17} $ denote $\bar{u}_{i-1,j-1},\bar{u}_{i,j-1},\bar{u}_{i+1,j-1},\bar{u}_{i-1,j},\bar{u}_{i,j},\bar{u}_{i+1,j},\bar{u}_{i-1,j+1},\bar{u}_{i,j+1},\bar{u}_{i+1,j+1},\\
	\bar{v}_{i-1,j},\bar{v}_{i+1,j},\bar{v}_{i,j-1},\bar{v}_{i,j+1},\bar{w}_{i-1,j},\bar{w}_{i+1,j},\bar{w}_{i,j-1}$,and $\bar{w}_{i,j+1}$, respectively.